\documentclass{article}

\usepackage{amsmath}
\usepackage{amsfonts}
\usepackage{amssymb}
\usepackage{amsthm}
\usepackage{graphics}
\usepackage{enumerate}
\usepackage{color}
\usepackage{tikz-cd}[2010/10/13] 
\usepackage[bookmarks=false]{hyperref}
\usepackage{cleveref}
\usepackage{tensor}
\usepackage{booktabs}
\usepackage[normalem]{ulem}
\usetikzlibrary{arrows}

\newcommand{\D}{\nabla}
\newcommand{\ed}{\mathrm{d}}

\newcommand{\R}{\mathbb{R}}

\newcommand{\ol}{\overline}

\newtheorem{corollary}{Corollary}
\newtheorem{proposition}{Proposition}
\newtheorem{theorem}{Theorem}
\newtheorem{lemma}{Lemma}
\newtheorem{remark}{Remark}
\theoremstyle{definition}
\newtheorem{definition}{Definition}
\newtheorem{example}{Example}[section]
\newtheorem*{remark*}{Remark}
\DeclareMathOperator{\Image}{Im}

\newcommand{\itoLTwoProjectedEquation}{\begin{equation*}\begin{split}
\ed\theta ^1 &= \left(-\frac{1}{4} \theta ^1 \left(\theta ^2\right)^2 \left(3 \epsilon ^2 \left(4 \left(\theta ^1\right)^4-4 \left(\theta ^2\right)^2 \left(\theta ^1\right)^2-3 \left(\theta ^2\right)^4\right)+16 \epsilon  \left(\theta ^1\right)^2+4\right)\right) \, \ed t \\
&{}+ \left(\frac{1}{2} \left(\theta ^2\right)^2 \left(3 \epsilon  \left(2 \left(\theta ^1\right)^2+\left(\theta ^2\right)^2\right)+2\right)\right) \, \ed Y_t \\ 
\ed\theta ^2 &= \left(-\frac{9 \epsilon ^2 \left(\theta ^2\right)^8+\left(\theta ^2\right)^4 \left(60 \epsilon ^2 \left(\theta ^1\right)^4+48 \epsilon  \left(\theta ^1\right)^2+4\right)+6 \epsilon  \left(\theta ^2\right)^6 \left(9 \epsilon  \left(\theta ^1\right)^2+2\right)-4}{8 \theta ^2}\right) \, \ed t \\
&{}+ \left(3 \epsilon  \theta ^1 \left(\theta ^2\right)^3\right) \, \ed Y_t
\end{split}\end{equation*}}
\newcommand{\itoHellingerProjectedEquation}{\begin{equation*}\begin{split}
\ed\theta ^1 &= \left(-\theta ^1 \left(\theta ^2\right)^2 \left(3 \epsilon ^2 \left(\left(\theta ^1\right)^4+4 \left(\theta ^2\right)^2 \left(\theta ^1\right)^2+6 \left(\theta ^2\right)^4\right)+\epsilon  \left(4 \left(\theta ^1\right)^2+6 \left(\theta ^2\right)^2\right)+1\right)\right) \, \ed t \\
&{}+ \left(\left(\theta ^2\right)^2 \left(3 \epsilon  \left(\left(\theta ^1\right)^2+\left(\theta ^2\right)^2\right)+1\right)\right) \, \ed Y_t \\ 
\ed\theta ^2 &= \left(-\frac{27 \epsilon ^2 \left(\theta ^2\right)^8+\left(\theta ^2\right)^4 \left(15 \epsilon ^2 \left(\theta ^1\right)^4+12 \epsilon  \left(\theta ^1\right)^2+1\right)+9 \epsilon  \left(\theta ^2\right)^6 \left(6 \epsilon  \left(\theta ^1\right)^2+1\right)-1}{2 \theta ^2}\right) \, \ed t \\
&{}+ \left(3 \epsilon  \theta ^1 \left(\theta ^2\right)^3\right) \, \ed Y_t
\end{split}\end{equation*}}
\newcommand{\itoADFEquation}{\begin{equation*}\begin{split}
\ed\theta ^1 &= \left(-\theta ^1 \left(\theta ^2\right)^2 \left(3 \epsilon ^2 \left(\left(\theta ^1\right)^4+4 \left(\theta ^2\right)^2 \left(\theta ^1\right)^2+3 \left(\theta ^2\right)^4\right)+\epsilon  \left(4 \left(\theta ^1\right)^2+6 \left(\theta ^2\right)^2\right)+1\right)\right) \, \ed t \\
&{}+ \left(\left(\theta ^2\right)^2 \left(3 \epsilon  \left(\left(\theta ^1\right)^2+\left(\theta ^2\right)^2\right)+1\right)\right) \, \ed Y_t \\ 
\ed\theta ^2 &= \left(-\frac{9 \epsilon ^2 \left(\theta ^2\right)^8+\left(\theta ^2\right)^4 \left(15 \epsilon ^2 \left(\theta ^1\right)^4+12 \epsilon  \left(\theta ^1\right)^2+1\right)+3 \epsilon  \left(\theta ^2\right)^6 \left(15 \epsilon  \left(\theta ^1\right)^2+2\right)-1}{2 \theta ^2}\right) \, \ed t \\
&{}+ \left(3 \epsilon  \theta ^1 \left(\theta ^2\right)^3\right) \, \ed Y_t
\end{split}\end{equation*}}

\newcommand{\stratLTwoProjectedEquation}{\begin{equation*}\begin{split}
\ed\theta ^1 &= \left(-\frac{1}{4} \theta ^1 \left(\theta ^2\right)^2 \left(3 \epsilon ^2 \left(4 \left(\theta ^1\right)^4-4 \left(\theta ^2\right)^2 \left(\theta ^1\right)^2-3 \left(\theta ^2\right)^4\right)+16 \epsilon  \left(\theta ^1\right)^2+4\right)\right) \, \ed t \\
&{}+ \left(\frac{1}{2} \left(\theta ^2\right)^2 \left(3 \epsilon  \left(2 \left(\theta ^1\right)^2+\left(\theta ^2\right)^2\right)+2\right)\right) \, \ed Y_t \\ 
\ed\theta ^2 &= \left(-\frac{47 \epsilon ^2 \left(\theta ^2\right)^8+\left(\theta ^2\right)^4 \left(60 \epsilon ^2 \left(\theta ^1\right)^4+48 \epsilon  \left(\theta ^1\right)^2+4\right)+2 \epsilon  \left(\theta ^2\right)^6 \left(33 \epsilon  \left(\theta ^1\right)^2+8\right)-4}{8 \theta ^2}\right) \, \ed t \\
&{}+ \left(3 \epsilon  \theta ^1 \left(\theta ^2\right)^3\right) \, \ed Y_t
\end{split}\end{equation*}}
\newcommand{\stratHellingerProjectedEquation}{\begin{equation*}\begin{split}
\ed\theta ^1 &= \left(-\theta ^1 \left(\theta ^2\right)^2 \left(3 \epsilon ^2 \left(\left(\theta ^1\right)^4+4 \left(\theta ^2\right)^2 \left(\theta ^1\right)^2+6 \left(\theta ^2\right)^4\right)+\epsilon  \left(4 \left(\theta ^1\right)^2+6 \left(\theta ^2\right)^2\right)+1\right)\right) \, \ed t \\
&{}+ \left(\left(\theta ^2\right)^2 \left(3 \epsilon  \left(\left(\theta ^1\right)^2+\left(\theta ^2\right)^2\right)+1\right)\right) \, \ed Y_t \\ 
\ed\theta ^2 &= \left(-\frac{36 \epsilon ^2 \left(\theta ^2\right)^8+\left(\theta ^2\right)^4 \left(15 \epsilon ^2 \left(\theta ^1\right)^4+12 \epsilon  \left(\theta ^1\right)^2+1\right)+9 \epsilon  \left(\theta ^2\right)^6 \left(6 \epsilon  \left(\theta ^1\right)^2+1\right)-1}{2 \theta ^2}\right) \, \ed t \\
&{}+ \left(3 \epsilon  \theta ^1 \left(\theta ^2\right)^3\right) \, \ed Y_t
\end{split}\end{equation*}}
\newcommand{\ekfEquation}{\begin{equation*}\begin{split}
\ed\theta ^1 &= \left(\left(\theta ^2\right)^2 \left(-\left(3 \epsilon  \left(\theta ^1\right)^2+1\right)\right) \left(\theta ^1+\epsilon  \left(\theta ^1\right)^3\right)\right) \, \ed t \\
&{}+ \left(\left(\theta ^2\right)^2 \left(3 \epsilon  \left(\theta ^1\right)^2+1\right)\right) \, \ed Y_t \\ 
\ed\theta ^2 &= \left(\frac{1-\left(\theta ^2\right)^4 \left(3 \epsilon  \left(\theta ^1\right)^2+1\right)^2}{2 \theta ^2}\right) \, \ed t \\
&{}+ \left(0\right) \, \ed Y_t
\end{split}\end{equation*}}
\newcommand{\itoLTwoProjectedEquationB}{\begin{equation*}\begin{split}
\ed\theta ^1 &= \left(-\frac{1}{4} \theta ^1 \left(\theta ^2\right)^2 \left(3 \epsilon ^2 \left(4 \left(\theta ^1\right)^4-4 \left(\theta ^2\right)^2 \left(\theta ^1\right)^2-9 \left(\theta ^2\right)^4\right)+16 \epsilon  \left(\theta ^1\right)^2+4\right)\right) \, \ed t \\
&{}+ \left(\frac{1}{2} \left(\theta ^2\right)^2 \left(3 \epsilon  \left(2 \left(\theta ^1\right)^2+\left(\theta ^2\right)^2\right)+2\right)\right) \, \ed Y_t \\ 
\ed\theta ^2 &= \left(\frac{3 \epsilon ^2 \left(\theta ^2\right)^8-4 \left(\theta ^2\right)^4 \left(15 \epsilon ^2 \left(\theta ^1\right)^4+12 \epsilon  \left(\theta ^1\right)^2+1\right)-2 \epsilon  \left(\theta ^2\right)^6 \left(15 \epsilon  \left(\theta ^1\right)^2+2\right)+4}{8 \theta ^2}\right) \, \ed t \\
&{}+ \left(3 \epsilon  \theta ^1 \left(\theta ^2\right)^3\right) \, \ed Y_t
\end{split}\end{equation*}}
\newcommand{\itoHellingerProjectedEquationB}{\begin{equation*}\begin{split}
\ed\theta ^1 &= \left(-\theta ^1 \left(\theta ^2\right)^2 \left(3 \epsilon ^2 \left(\left(\theta ^1\right)^4+4 \left(\theta ^2\right)^2 \left(\theta ^1\right)^2+3 \left(\theta ^2\right)^4\right)+\epsilon  \left(4 \left(\theta ^1\right)^2+6 \left(\theta ^2\right)^2\right)+1\right)\right) \, \ed t \\
&{}+ \left(\left(\theta ^2\right)^2 \left(3 \epsilon  \left(\left(\theta ^1\right)^2+\left(\theta ^2\right)^2\right)+1\right)\right) \, \ed Y_t \\ 
\ed\theta ^2 &= \left(-\frac{18 \epsilon ^2 \left(\theta ^2\right)^8+\left(\theta ^2\right)^4 \left(15 \epsilon ^2 \left(\theta ^1\right)^4+12 \epsilon  \left(\theta ^1\right)^2+1\right)+3 \epsilon  \left(\theta ^2\right)^6 \left(15 \epsilon  \left(\theta ^1\right)^2+2\right)-1}{2 \theta ^2}\right) \, \ed t \\
&{}+ \left(3 \epsilon  \theta ^1 \left(\theta ^2\right)^3\right) \, \ed Y_t
\end{split}\end{equation*}}

\title{Optimal approximation of SDEs on submanifolds:\\ the It\^o-vector and It\^o-jet projections}
\author{John Armstrong \\ Dept.\ of Mathematics \\ King's College London \\ {\small \tt{john.1.armstrong@kcl.ac.uk}} \and Damiano Brigo \\ Dept.\ of Mathematics \\ Imperial College London \\ {\small \tt{damiano.brigo@imperial.ac.uk}}}

\date{First version: Sept. 12, 2015. This version: \today}

\begin{document}

\maketitle
\begin{abstract}

We define two new notions of projection of a stochastic differential equation (SDE) onto a submanifold: the It\^o-vector and It\^o-jet projections. This allows one to systematically develop low dimensional approximations to high dimensional SDEs using differential geometric techniques. The approach generalizes the notion of projecting a vector field onto a submanifold in order to derive approximations to ordinary differential equations, and improves the previous Stratonovich projection method by adding  optimality analysis and results. Indeed, just as in the case of ordinary projection, our definitions of projection are based on optimality arguments and give in a well-defined sense ``optimal'' approximations to the original SDE in the mean-square sense. We also show that the Stratonovich projection satisfies an optimality criterion that is more ad hoc and less appealing than the criteria satisfied by the It\^o projections we introduce.

As an application we consider approximating the solution of the non-linear filtering problem with a Gaussian distribution
and show how the newly introduced It\^o projections lead to optimal approximations in the Gaussian family and briefly discuss the optimal approximation for more general families of distribution. We perform a numerical comparison of
our optimally approximated filter with the classical Extended Kalman Filter to demonstrate the efficacy of the approach.
\end{abstract}

\bigskip

{\bf Keywords:} Stochastic differential equations, Jets, SDEs as jets, SDEs projection on a manifold, SDEs on submanifolds, Stratonovich projection, It\^o-vector projection, It\^o-jet projection,  Optimal projection, Gaussian It\^o-jet filter, Gaussian It\^o-vector filter.  

\bigskip

{AMS classification codes: 58A20, 39A50, 58J65, 60H10, 60J60, 65D18}

\bigskip

\section{Introduction}

In this paper we define three notions of projecting a stochastic differential equation (SDE)
onto a (sub)manifold $M$. Our aim is to derive practical numerical methods for solving SDEs and we will illustrate our theory with an example drawn from signal processing. 

\begin{figure}[h!]
\center{\includegraphics[scale=0.35]{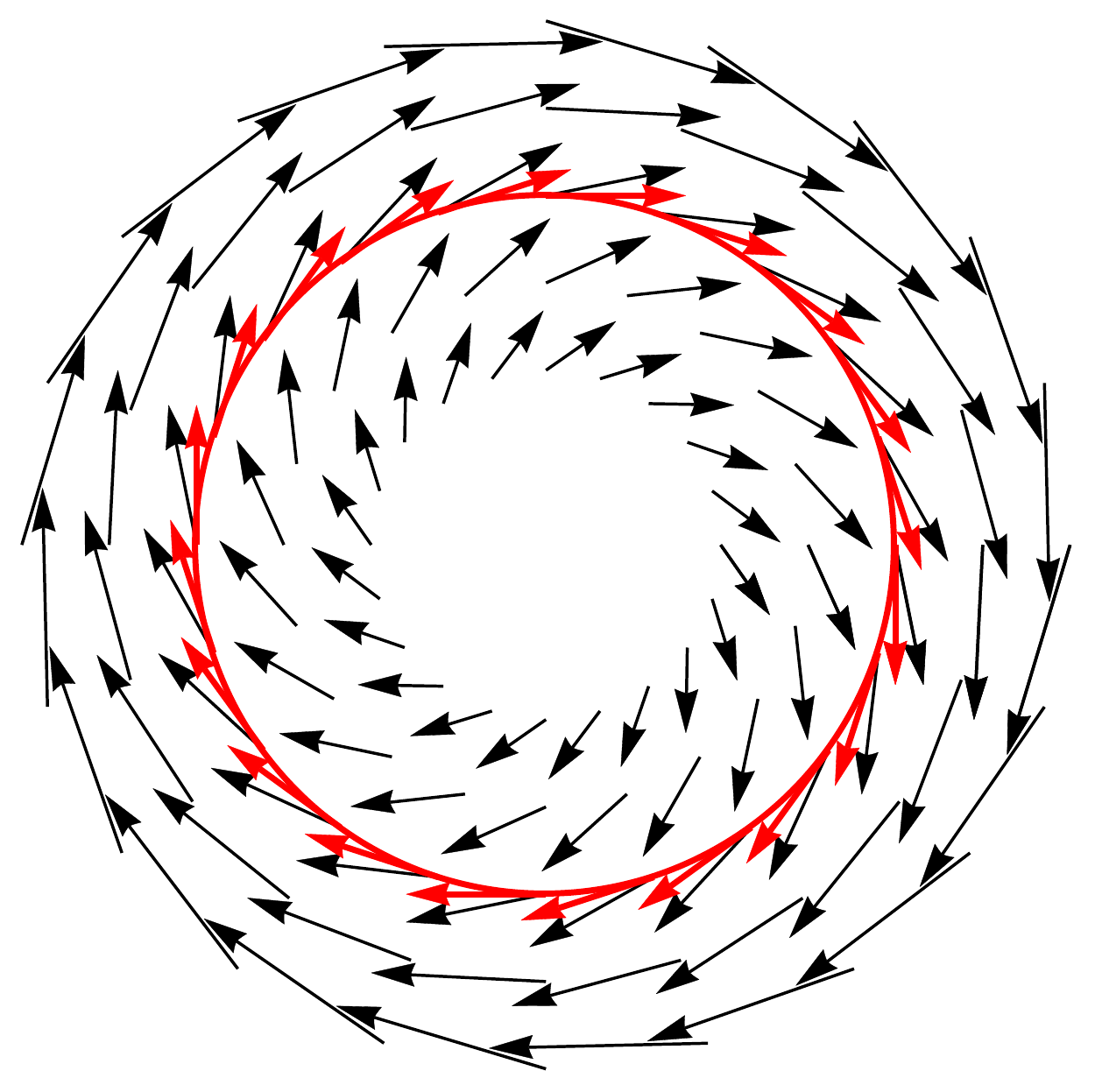}\hspace{1cm} \includegraphics[scale=0.35]{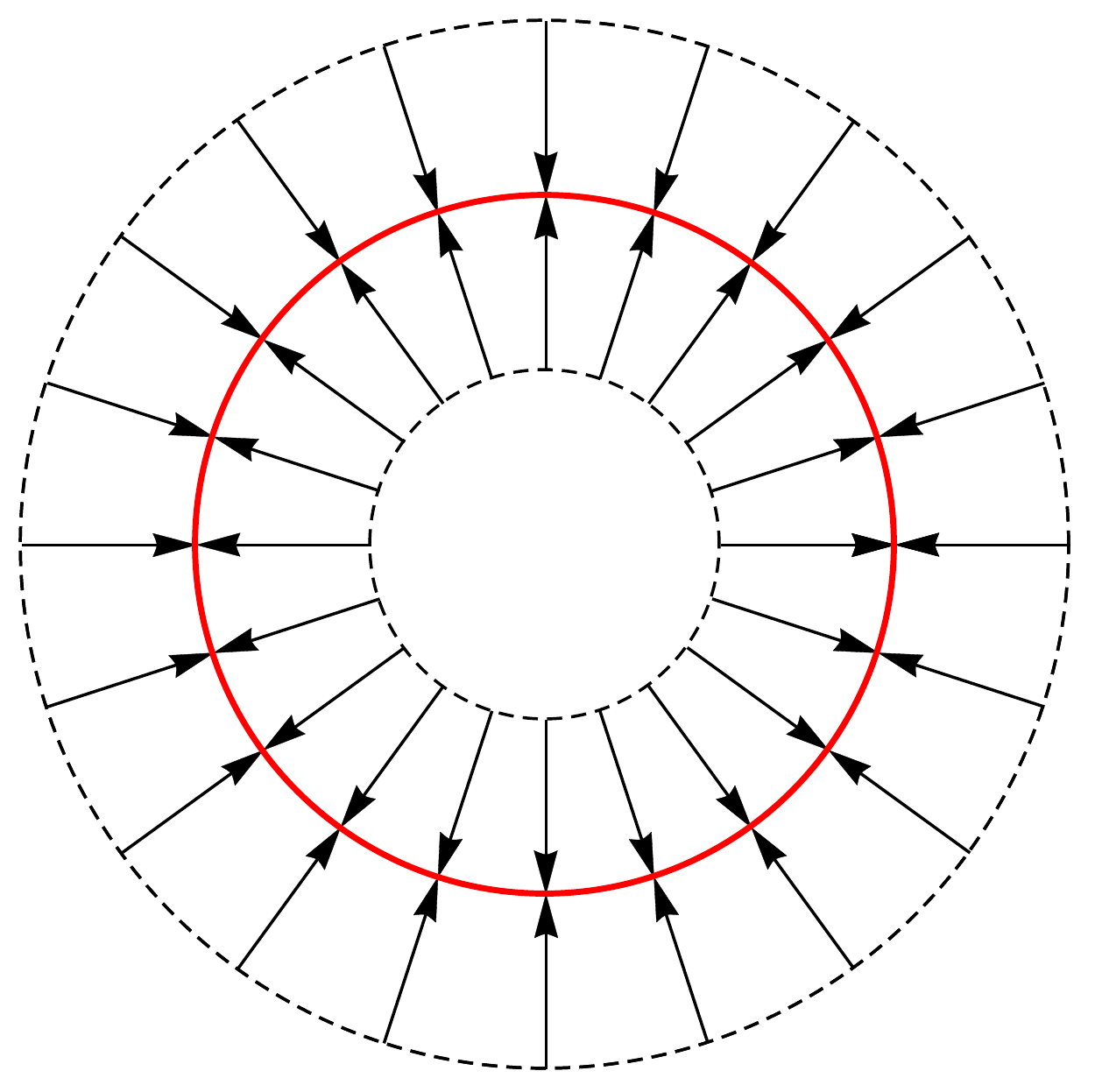}}
\caption{Left: A pictorial representation of the projection of an ODE defined on $\R^2$ to an ODE defined on the circle. Right: Metric projection from $\R^2$ to the circle.}
\label{figure:odeProjection}
\label{figure:metricProjection}
\end{figure}

To explain the general idea, let us first consider projecting an ordinary differential equation (ODE) from the Euclidean space $\R^r$ onto an $n$-dimensional manifold $M \subseteq \R^r$.
An ODE in $\R^r$ can be thought of as defining a vector field in $\R^r$. At every point $x \in M$ we
can use the Euclidean metric to project the vector at $x$ onto the tangent space $T_xM$. In this
way one obtains a vector field on $M$ which can be thought of as a new ODE on $M$ that
approximates the full ODE in $\R^r$. This is illustrated in \Cref{figure:odeProjection}. It is easy to prove that this will be the
best way of approximating the ODE in $\R^r$ with an ODE on $M$. To be precise, if the initial condition for an ODE is a point $x$ on the manifold, then any curve on $M$ with tangent not equal to the projected vector field will diverge from the solution to the ODE faster than a curve which is tangent to the projected vector field. In this sense, the projected ODE is the only ODE which is ``optimal'' at each point.
This paper addresses the question of how projection can be generalized from ODEs to SDEs. After some brief preliminaries on It\^o--Taylor series in \Cref{section:taylorSeries}, we answer this question by describing three possible generalizations to SDEs in \Cref{section:projectionSection}.

The first generalization of projection to SDEs has been proposed previously: what we shall call the Stratonovich projection.
The Stratonovich projection is obtained by simply applying the projection operator to the coefficients of the SDE written in Stratonovich form. No optimality result has been derived for the Stratonovich projection. This projection has simply been derived heuristically from the deterministic case. Nevertheless, it appears to be a good approximation in practice and it has been used to find good quality numerical solutions to the non-linear filtering problem (See \cite{brigo1}, \cite{brigo2}, \cite{armstrongBrigo}). The Stratonovich projection is a natural first choice from the following point of view. As is obvious to anyone with experience of stochastic differential equations on manifolds, and we refer to the monographs and articles \cite{emery},\cite{elworthy},\cite{gliklikh}, \cite{brzezniak}, \cite{hsu}, \cite{armstrongBrigoJets}, simply applying the projection operator to the coefficients of the SDE written in It\^o form will not work. This is because solutions to the projected equation don't stay on the manifold, contrary to the Stratonovich case. Nevertheless, we will be able to obtain two modifications of this idea, which we will call the {\em It\^o-vector projection} and the {\em It\^o-jet projection}. These both give well-defined SDEs on the manifold.

We derive the It\^o-vector projection by seeking an SDE on the manifold which
optimally approximates the original SDE on the manifold when the size of the
errors are measured in the mean square ambient metric of $\R^r$. The mean squared error between a trajectory following the original SDE and following an SDE on the manifold will typically grow at a rate $O(t^{\frac{1}{2}})$. The diffusion term of the projected SDE is determined by minimizing the coefficient of $t^\frac{1}{2}$ in this growth estimate. Choosing the drift term is more delicate, but we give two minimization arguments that indicate that the optimal choice of drift term is given by what we call the It\^o-vector projection.
The first argument identifies the drift by minimizing the coefficient of the $O(t)$ term in the estimate of the error, notwithstanding the fact that there is also an $O(t^{\frac{1}{2}})$ term. The second argument is to find an SDE on the manifold such that the difference between the means of the solutions to the original SDE and to the SDE on the manifold are minimized.

Both of these arguments are somewhat unsatisfying. As an alternative approach we consider finding the SDE on the manifold that most closely tracks the {\em metric projection} of the solution to the original SDE. The metric projection is the map that sends a point in the ambient space to the closest point of an embedded manifold $M$. It is well known to be well-defined and smooth on a tubular neighbourhood. The metric projection is illustrated in \Cref{figure:metricProjection}.
It is possible to find an SDE on the manifold such that the mean squared distance between the solutions on the manifold and the metric projection of the solution to the original SDE grows at a rate $O(t)$. This requirement determines the diffusion term of the SDE on the manifold and makes the $O(t^{\frac{1}{2}})$ term coefficient vanish, rather than merely minimize it. Minimizing the coefficient of the order $t$ term in this estimate determines the drift.
We call the SDE determined in this way the {\em It\^o-jet projection}.

It is natural to ask if the Stratonovich projection can also be derived from an optimality argument. We will show that the Stratonovich projection is optimal when using a time-reflection-symmetric optimality criterion anchored to the deterministic intial condition of the process as a special state. We will clarify this notion of optimality in the paper.
 However, as we will see, for our applications to filtering, the form of optimality achieved by the Stratonovich projection is not particularly useful. This is because the filtering problem is inherently asymmetric in time, as indeed are most applications of SDEs. Nevertheless, it is conceivable that in some applications of SDEs to physics, time reversal symmetry may be a paramount concern. In this case the Stratonovich projection may be preferred.

Surprisingly the It\^o-vector projection, the It\^o-jet projection and the Strat\-ono\-vich projection are all distinct. All of them reduce to classical projection in the case of ODEs. Thus, while optimality arguments lead to a single best method for projecting ODEs, the situation is more complex for SDEs. Since both It\^o projections are derived from optimality arguments that are much less ad hoc than the argument for optimality of the Stratonovich projection, there is a clear sense in which they are an improvement upon the Stratonovich projection---both theoretically and in practice.

However, it is not immediately
clear whether one should prefer the It\^o-vector or the It\^o-jet projection.
We investigate this question in \Cref{section:lowDimensionalExample}.

We consider a simple toy example in 
\Cref{section:lowDimensionalExample} which we believe strongly suggests
that the It\^o-jet projection is the better approximation. We also prove a simple theorem that shows how this example can be generalized.

We use this same toy example to illustrate another (entirely non-rigorous) reason for preferring the It\^o-jet projection: mathematical aesthetics. As we shall see, each of the different notions of projection is 
best understood using different formulations of SDEs on manifolds. As its name suggests,
the Stratonovich projection is most readily understood using Stratonovich calculus. The It\^o-vector
projection is most readily understood using the formulation of SDEs on manifolds in terms of It\^o
calculus first introduced by It\^o in \cite{ito2}. Finally, the It\^o-jet projection is most readily understood using the $2$-jet formulation of \cite{armstrongBrigoJets}. As we will see,
the It\^o-jet projection has a very elegant formulation in the language of $2$-jets. It is even possible to draw a diagram that allows one to interpret the It\^o-jet projection visually. We will present a diagram that visually represents the It\^o-jet projection of our toy example.
In fact, the
development of the $2$-jet formulation of SDEs in \cite{armstrongBrigoJets} was originally motivated by
the development of these projection methods. It is for this reason that we have called the projections the It\^o-vector and It\^o-jet projections respectively.

\Cref{section:secondOrderProjection} is devoted to a detailed calculation
of the It\^o-jet projection in local coordinates. This calculation amounts
to computing the Taylor series for the metric projection map up to second
order. This calculation is essential to using the projection for applications.

\Cref{applicationSection} demonstrates how the notion of projection can be applied in 
practice.  In particular, we will apply it to the non-linear filtering problem.  We will derive general projection
formulae for the non-linear filtering problem. We will then apply this to the problem of approximating a non-linear filter using a  Gaussian distribution. A reader who is unfamiliar with non-linear filtering will want to consult  \Cref{filteringSubsection} for a brief review.

Gaussian approximations to non-linear filters are widely used in practice (see for example \cite{jazwinski,crisan}). In particular,  the Extended Kalman Filter (EKF) is a popular approximation technique. Other Gaussian approximations exist such as Assumed Density Filters (ADF) and filters derived from the Stratonovich projection. Our theory indicates that all these classical techniques can be improved upon by using the It\^o projections (at least over small time intervals). We confirm this with a numerical example.

The utility of the projection method is by no means restricted to the
filtering problem nor to such simple approximations as Gaussian filters. Our previous work \cite{armstrongbrigomcss} shows how the Stratonovich projection can be used to generate far more sophisticated filters and it is clear that the idea of projection should be widely applicable in the study of ODEs, SDEs,  PDEs and SPDEs.
Nevertheless by focussing on Gaussian filters we can examine in detail the idea that there may be many useful ways of approximating an SDE on a submanifold, but that the It\^o projections are in some sense optimal
amongst these approximation methods. The point we wish to emphasize is that the It\^o projections are able to tell us something new even about the well-worn topic of approximating the non-linear filtering problem using Gaussian distributions.

Finally in \Cref{section:conclusions} we summarize our findings.

\section{Stochastic Taylor Series}
\label{section:taylorSeries}
The main technical tool we will use are stochastic Taylor series.
These are described in detail in \cite{kloedenAndPlaten}. In this
section we will recall the main definitions and results. We
will make some minor notational changes so that we can use the Einstein
summation convention.

Let $X_t$ satisfy a $d$-dimensional stochastic differential equation
driven by $m$ independent Brownian motions $W^\alpha_t$, $\alpha=1,2,\ldots,m$. 
We write
\begin{equation}
\ed X_t = a(X,t) \ed t + b_\alpha(X,t) \ed W^\alpha_t
\label{generalSDE}
\end{equation}
where $X_t$ is an random process taking values on $\R^d$. $a$ and
$b_\alpha$ are also $\R^d$ valued for each $\alpha$.
We are using the Einstein
summation convention that when there are matching indices
in an expression one should take the sum over the given index. Thus
\eqref{generalSDE} is an abbreviation for:
\begin{equation*}
\ed X_t = a(X,t) \ed t + \sum_{\alpha=1}^m b_\alpha(X,t) \ed W^\alpha_t.
\end{equation*}
The advantage of the Einstein summation convention is not simply
that it makes formulae shorter. The convention also makes it easier to spot
incorrect formulae. This is because, in formulae that are valid in
all coordinate systems, the summed indices should always consist of
one upper and one lower index.

In this section we will use Greek indices to index the different Brownian motions and Roman indices to index components of vectors in $\R^d$. This additional convention is not strictly necessary as the range of the index can
be deduced from the position of the index alone.

A multi-index $\xi$ is defined to be a finite list of integer numbers between $0$ and $m$ and this definition includes the empty list $()$. Let $l(\xi)$ denote the length of $\xi$. Let $n(\xi)$ denote the
number of zeros in $\xi$. For $\xi$ with length greater than 0, we define:
$-\xi$ to be the result of removing the first element from $\xi$; $\xi-$
for the result of removing the last element;  $\xi_1$ for the first element; 
and $\xi_{-1}$ for the last element.

Multi-indices enumerate stochastic integrals with respect to the
Brownian motions $W^\alpha_t$ and time. The following definitions are related to those on page 169 of \cite{kloedenAndPlaten}. We define $W^0_t:=t$ so that
the indices equal to 0 correspond to time.
We define the multi-integral associated with $\xi$ by:
\[
I_{t_1,t_2}^\xi(f) =
\begin{cases}
f(t_2) & \hbox{if } l(\xi)=0 \\
\int_{t_1}^{t_2} I_{t_1,s}^{\xi-} \ed W^{\xi_{-1}}_s & \hbox{otherwise}.
\end{cases}
\]
For example the multi-index
$(0,1,2)$ is associated with integrating with respect first to time,
then $W^1$, then $W^2_t$.
\[
I_{t_1,t_2}^{(0,1,2)} (f) = \int_{t_1}^{t_2} \int_{t_1}^{u} \int_{t_1}^{v}
f(w) \, \ed w \, \ed W^1_{v} \, \ed W^2_{u}.
\]
We re-express the notation in \cite{kloedenAndPlaten} (page 177, Eqs. 3.1--3.3) by  defining differential operators $L_\xi$ associated to a multi index as:
\[
L_\xi f = 
\begin{cases}
f & \hbox{if } l(\xi)=0 \\
\frac{\partial f}{\partial t} + a^i \frac{\partial f}{\partial x^i}
+ \frac{1}{2} b^i_\alpha b^j_\beta g^{\alpha \beta}_E \frac{ \partial^2
f}{\partial x^i \partial x^j } & \hbox{if } \xi=(0) \\
b^i_{\xi_1} \frac{\partial f}{\partial x^i} & \hbox{if } l(\xi)=1 \hbox{ and } \xi\neq(0) \\
L_{\xi_1} (L_{-\xi} f) & \hbox{otherwise}.
\end{cases}
\]
Here $g^{\alpha \beta}_E$ denotes the covariance matrix of the $d$ Brownian
motions $W^\alpha_t$. Since we have assumed that the Brownian motions are independent, this will equal the identity matrix. We choose to write
$g_E$ instead of using the Kronecker delta because it transforms as a tensor of type $(2,0)$. In addition, one can simply replace $g_E$ with
the quadratic co-variation tensor if one wishes to consider
SDEs driven by more general continuous semi-martingales.

Since $L_\xi$ contains a total of $l(\xi)+n(\xi)$ derivatives,
$L_\xi$ acts on functions in $C^{l(\xi)+n(\xi)}( \R^d \times  \R^+, \R )$.

The following definition is related to Eq. (9.1) page 206 in \cite{kloedenAndPlaten}.
\begin{definition}
The {\em It\^o--Taylor expansion} of order $\gamma=0,\frac{1}{2},1,\ldots$
is given by
\[
X^\gamma_t = \sum_{l(\xi)+n(\xi) \leq 2 \gamma} L_\xi(x)\rvert_{(t_0,X_{t_0})}
I^\xi_{t_0,t}(1)
\]
where $x$ denotes the function $x(t,X)=X$. When we speak of the
expansion of a given order, we will assume that all the necessary
derivatives exist.
\end{definition}
The It\^o--Taylor expansion allows one to approximate $X_t$ using $X^\gamma_t$.
Loosely speaking, this approximation will be accurate in mean
squared up to order $\gamma$. A precise statement is given in \Cref{proposition:strongConvergence}.
\begin{definition}
The {\em weak It\^o--Taylor expansion} of order $\beta=0,1,2,\ldots$
is given by
\[
\eta_\beta(t) = \sum_{l(\xi) \leq \beta} L_\xi(x)\rvert_{(t_0,X_{t_0})}
I^\xi_{t_0,t}(1)
\]
where $x$ denotes the function $x(t,X)=X$. When we speak of the
expansion of a given order, we will assume that all the necessary
derivatives exist.
\end{definition}
The weak It\^o--Taylor expansion is of interest if one measures
the error using the size of the expectation of the error, rather than
the expectation of the size of the error. We will give a precise
statement in \Cref{proposition:weakConvergence}.

Given a smooth vector valued function $f$ defined on $\R^d$ we have by It\^o's lemma that
\begin{equation}
\ed f(X_t) = L_0(f) \ed t + L_\alpha(f) \ed W^\alpha_t.
\label{fSDE}
\end{equation}
The system of equations \eqref{generalSDE} and \eqref{fSDE} define a higher
dimensional SDE. We can use this to compute It\^o--Taylor expansions for this
higher dimensional system and hence compute approximations to $f(X_t)$. 
This calculation gives rise to the following more general definition.
\begin{definition}
The {\em It\^o--Taylor expansion} of order $\gamma=0,\frac{1}{2},1,\ldots$
for $f(X_t)$ is given by
\[
f^\gamma_t = \sum_{l(\xi)+n(\xi) \leq 2 \gamma} L_\xi(f)\rvert_{(t_0,X_{t_0})}
I^\xi_{t_0,t}(1)
\]
When we speak of the
expansion of a given order, we will assume that all the necessary
derivatives exist. The weak It\^o--Taylor expansion for $f(X_t)$
is defined similarly.
\end{definition}

\begin{lemma}
\label{momentCalculations}
We suppose that for all $i$, $W^i_0=0$.
Given a time $t$, and $i$, $j \in \{1, \ldots m\}$, the integrals
\[
\begin{split}
I^{(0)}_{0,t}(1) &= t \\
I^{(i)}_{0,t}(1) &= W^i_t \\
I^{(i,j)}_{0,t}(1) &= \int_0^t W^i_s \ed W^j_s \\
\end{split}
\]
are orthogonal in expectation.
\end{lemma}
\begin{proof}
We first show that
\begin{equation}
E\left( I^{(i)}_{0,t}(1) I^{(j,k)}_{0,t}(1)\right)
= E\left( W^i_t \int_0^t W^j_s \ed W^k_s \right) = 0
\label{orthogonalityRelation}
\end{equation}
for all $i$, $j$, $k \in \{1, \ldots m\}$.
If $i\neq j$ then we see, by reversing the sign of $W^j$, that 
\[ E\left( W^i_t \int_0^t W^j_s \ed W^k_s \right) =
-E\left(W^i_t \int_0^t W^j_s \ed W^k_s \right)
\]
Hence \eqref{orthogonalityRelation} is zero unless $i=j$. The same
argument shows  \eqref{orthogonalityRelation}  is zero unless $i=k$.
Finally we note that when $i=j=k$, \eqref{orthogonalityRelation} simplifies
to
\[
E\left( W^i_t \int W^i_s \ed W^i_s \right) = E\left(\int \ed W^i_s \int W^i_s \ed W^i_s \right)
= E\left(\int_0^t W^i_s \ed s\right) = 0 
\]
by the It\^o isometry.

We also need to show that if $i\neq k$ or $j\neq l$
\begin{equation*}
E\left( I^{(i,j)}_{0,t}(1) I^{(k,l)}_{0,t}(1) \right) = 0.
\end{equation*}
This follows from Lemma 5.7.2 on page 191 of \cite{kloedenAndPlaten}.

The other cases are trivial.
\end{proof}

For completeness, we wish to state some results on the convergence
of It\^o--Taylor series. We will first
need a few more definitions.

First we define spaces ${\cal H}_\xi$ associated with
multi-indices $\xi$. Associated to the empty index $()$ we
have the set ${\cal H}_{()}$ of adaptad cadlag processes $f_t$ with
\[ |f(t,\omega)|<\infty \]
with probability one for each $t\geq 0$. ${\cal H}_{(0)}$
consists of the adapted cadlag processes with
\[
\int_0^t |f(s,\omega)| \ed s < \infty
\] 
with probability one for each $t\geq 0$. ${\cal H}_{(\alpha)}$
has the same definition for any positive $\alpha$: it is the
set of adapted cadlag processes with
\[
\int_0^t |f(s,\omega)|^2 \ed s < \infty
\] 
with probability one for each $t\geq 0$. We now recursively define
${\cal H}_\xi$ for $\xi$ of length greater than $1$ to be 
the set of adapted cadlag processes such that the integral
process $I^{\xi-}_{0,t}(f)$ when viewed
as a function of $t$ lies in ${\cal H}_{\xi_1}$.

We define ${\cal M}$ to be the set of all multi-indices.

Given a subset ${\cal A} \subseteq {\cal M}$ we define
the {\em remainder set} ${\cal B}({\cal A})$ to be the set
\[
{\cal B}({\cal A} ) = \{ \xi \in {\cal M} \setminus {\cal A} : -\xi \in {\cal A}.
\}
\]
Thus the remainder set contains all the indices immediately following
the indices in ${\cal A}$. By estimating integrals in the remainder
set, one can bound the error of the It\^o--Taylor series as we will
see below.

We define
\[ \Lambda_k = \{ \xi \in {\cal M}: l(\xi) + n(\xi) \leq k \}. \]
Thus the order $\gamma$ It\^o--Taylor series is a sum over multi-indexes
in $\Lambda_{2 \gamma}$.

We can now state a result on the convergence of the It\^o--Taylor
series. The following result is a simplified version of Proposition 5.9.1
in \cite{kloedenAndPlaten}.
\begin{proposition}
\label{proposition:strongConvergence}
Suppose that $L_{\xi} x\rvert_{X_{t_0},t_0} \in {\cal H}_\xi$ for all $\xi \in \Lambda_k$. Suppose that $L_{\xi} x\rvert_{X_{t},t} \in {\cal H}_\xi$
with
\[
\sup_{0\leq t \leq T} E\left(\left|(L_\xi x\rvert_{X_{t},t})\right|^2\right) \leq C_1
\]
for all $\xi \in {\cal B}(\Lambda_k)$ and some constant $C_1$. Then
\[ E\left( |X_t - X^{\frac{k}{2}}_t|^2 \right) \leq C_2 (t-t_0)^{k+1} \]
for some constant $C_2$. Here $X^{\frac{k}{2}}$ is the order $\frac{k}{2}$
It\^o-Taylor expansion with $k=0,1,\ldots$.
\end{proposition}

This next result on the convergence of weak It\^o--Taylor series
is a restatement of Proposition 5.11.1 in \cite{kloedenAndPlaten}.
\begin{proposition}
\label{proposition:weakConvergence}
Let $\beta \in \{1, 2, \ldots \}$ and $T \in (0, \infty)$ be given.
Let $C^l_P(\R^d, \R)$ denote the space of $l$ times continuously
differentiable functions whose derivatives of order up to and including
$l$ have polynomial growth. 
Suppose that $a^k$ and $b^{kj}$ are time-independent and satisfy Lipschitz conditions, linear growth bounds and belong to $C^{2(\beta+1)}_P(\R^d,\R)$.
Then for each $g \in C^{2(\beta+1)}_P(\R^d,\R)$ there exist constants
$K \in (0, \infty)$ and $r \in \{1, 2, \ldots \}$ such that
\[
\sup_{0\leq t \leq T} \left|
E \left(
g(X_t)-g(\eta_\beta(t)) 
\right) \right|
\leq K\left( 1 + |X_0|^{2 r}\right) T^{\beta+1}
\]
where $\eta_b(t)$ is the weak It\^o--Taylor series and the expectation
is taken conditional on the information at time $0$.
\end{proposition}

A quick note on how we plan to use the above results to obtain optimal approximations is in order. Take the strong Taylor series to make the point. In our applications the $X$ we will expand in strong Taylor series will be the difference between the true solution of a SDE and its approximation on a submanifold. Knowing that the error in the Taylor series is bounded as per Propositions \ref{proposition:strongConvergence} (and \ref{proposition:weakConvergence} for the weak case), we will then concentrate on minimizing the mean square of the truncated Taylor series.  Minimizing the mean square for the truncated expansion of the difference rather than the mean square of the difference itself will work in view of the convergence guaranteed by the above proposition. In this sense when we talk about ``minizing the error'' or ``difference'' later on in the paper we always mean minimizing the truncated expansion mean square.

\section{Projecting stochastic differential equations}
\label{section:projectionSection}

\begin{figure}[h!]
\center{\includegraphics{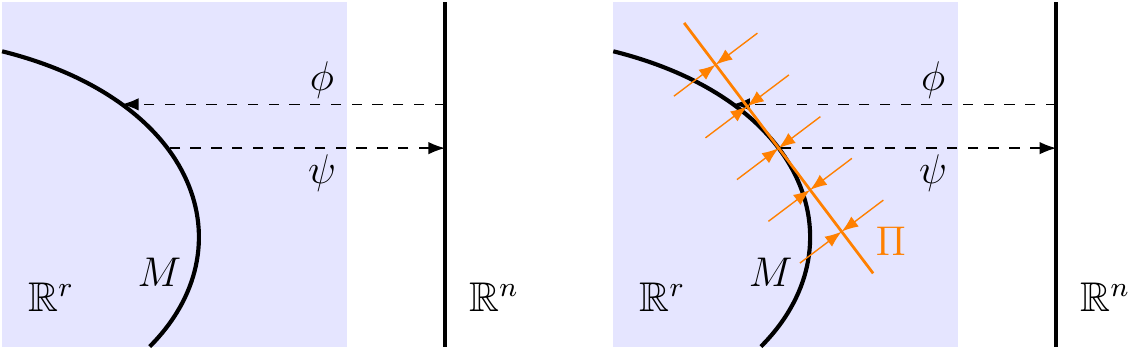}}
\caption{Left: An $n$-dimensional manifold $M$ in $\R^r$, $r>n$. Right: Tangent space linear projection used in the Stratonovich and It\^o-vector projections} 
\label{fig:phiAndPsi}
\end{figure}

Let $M$ be an $n$-dimensional submanifold of $\R^r$ with chart $\psi:U \to \R^n$ for some open neighbourhood $U$ in $M$. The inverse $\phi=\psi^{-1}$ gives
an embedding of $\Image \psi$ into $\R^r$. The setup is illustrated in \Cref{fig:phiAndPsi}.

Suppose we are given an It\^o SDE on $\R^r$, 
$\ed X_t = a(X_t,t) \, \ed t + b_\alpha(X_t,t) \, \ed W^\alpha_t$, that we write in concise form as

\begin{equation}
\ed X = a \, \ed t + b_\alpha \, \ed W^\alpha_t, \qquad X_0
\label{equationOnRr}
\end{equation}
with $X_0 \in M$.

We wish to find an SDE on $\R^n$ of the form $\ed Y_t = A(Y_t,t) \, \ed t + B_\alpha(Y_t,t) \, \ed W^\alpha_t$, again written concisely as
\begin{equation}
\ed Y = A \, \ed t + B_\alpha \ed W^\alpha_t, \qquad Y_0=\psi^{-1}(X_0),
\label{equationOnManifold}
\end{equation}
whose mapped solution $\phi(Y)$ in some sense approximates the solution $X$ of the original equation on $\R^r$.
We will consider three approaches.

\subsection{Stratonovich Projection}

\begin{definition}
Let $W_t$ be an ${\R^m}$ valued Brownian motion. Given a Stratonovich SDE on $\R^r$
\[ \ed X = \ol{a} \, \ed t + b_\alpha \, \circ \ed W^\alpha_t \]
and a chart $\psi:U \to \R^n$ for
some neighbourhood in $M$ we define the {\em Stratonovich projection} of
the SDE to be:
\[ \ed Y= {\ol A} \, \ed t + B_\alpha \circ \ed W^\alpha_t \]
where:
\begin{eqnarray}
{\ol A}(Y_t,t) &=& (\psi_*) \Pi_{\phi(Y_t)}( {\ol a}(\phi(Y_t),t)) \label{stratProjectionFulltermA} \\
B_\alpha(Y_t,t) &=& (\psi_*) \Pi_{\phi(Y_t)}( b_\alpha(\phi(Y_t),t))
\label{stratProjectionFull}
\end{eqnarray}
where $\Pi$ is the projection of $\R^r$ onto $\phi_*(\R^n)$ defined by the Euclidean metric.
\end{definition}
 
Because we know that projection of vector fields can be defined similarly, and because we know that the coefficients of Stratonovich SDEs transform like vector fields, we see that the definition above defines a Stratonovich SDE on $M$. Indeed, if one is willing to accept that projection of vector fields onto a submanifold is well-defined, then one could define the projection of a Stratonovich SDE as the projection of the coefficient functions.

Trying the same method for an It\^o SDE does not work. One cannot simply apply projection to the coefficient functions of an It\^o SDE because the coefficients of an It\^o SDE on a manifold do not transform like vector fields.

The {\em Stratonovich projection} of an It\^o SDE is trivially defined by the recipe:
\begin{enumerate}[(i)]
\item rewrite the It\^o SDE as a Stratonovich SDE;
\item apply the Stratonovich projection as defined above;
\item rewrite the resulting Stratonovich SDE as an It\^o SDE.
\end{enumerate}
In other words, while the definition of Stratonovich projection is
most conveniently expressed using Stratonovich calculus, the notion
of projection is independent of the calculus used to write down
the differential equations.

Linear projection provides the best possible way to approximate
vectors in $\R^r$ with vectors in $T_X M$. For ODEs, this implies that the projected ODE is the best possible approximation in $M$ of the original ODE. However, the situation is different for SDEs. It is not immediately clear how good an approximation the projected Stratonovich SDE solution $\phi(Y)$ is for the original SDE $X$ solution. For example, we cannot immediately extend the optimality argument for ODEs to Stratonovich SDEs pathwise, because of the rough paths property of SDEs solutions. In this sense, with the information we have given so far, the definition of the Stratonovich projection is motivated by purely
heuristic considerations.  Neverthless, the Stratonovich projection gives good results when applied to approximation of non-linear filtering problems (see \cite{brigo1}, \cite{brigo2}, \cite{armstrongBrigo}) and we will discuss optimality arguments later on, when discussing the It\^o-vector projection, and illustrate the time-symmetric optimality of the Stratonovich projection in detail. 

In the next sections we will use optimality arguments to derive two
alternative notions of projection.

\subsection{It\^o-vector projection}

We wish to consider the minimization problem of finding
coefficients $A$ and $B$ such that the solution of the SDE
\eqref{equationOnManifold}
has the property that $\phi(Y_t)$ is, in some sense, as close to 
the solution $X_t$ of \eqref{equationOnRr} as possible.

The next proposition shows how to give a precise meaning
to this notion using the It\^o--Taylor expansion.

\begin{proposition}
\label{generalProposition}
Let $f:\R^{d_x} \to \R^{d}$ and $F:\R^{d_y} \to \R^{d}$ be smooth maps. Let $x$ be a process on $\R^{d_x}$ and $y$ be a process on $\R^{d_y}$ given by:
\begin{equation}
\begin{split} 
\ed x_t &= a(x_t,t) \, \ed t + b_\alpha(x_t,t) \, \ed W^\alpha_t,\qquad x_0  \\
\ed y_t &= A(y_t,t) \, \ed t + B_\alpha(y_t,t) \, \ed W^\alpha_t,\qquad y_0  \\
\end{split}
\label{xyEquations}
\end{equation}
with $f(x_0)=F(y_0)$. Define
\[ z_t=f(x_t)-F(y_t). \]
Let $z^i_t$ denote 
the components of the order $i$ It\^o--Taylor expansion
for $z$. We have that:
\begin{equation}
\begin{split}
E(|z^{\frac{1}{2}}_t|^2)
&= \sum_{\alpha}|f_*(b_\alpha(x_0,0))- F_*(B_\alpha(y_0,0))|^2 t \\
E(|z^1_t|^2)
&= \sum_{\alpha}|f_*(b_\alpha(x_0,0))- F_*(B_\alpha(y_0,0))|^2 t \\
&\quad + \Bigg ( \Big|f_*(a(x_0,0)) - F_*(A(y_0,0)) \\
&\quad + \frac{1}{2} (\D_{b_{\alpha(x_0,0)}}f_*)b_\beta(x_0,0) g^{\alpha\beta}_E
 - \frac{1}{2} (\D_{B_{\alpha(y_0,0)}}F_*)B_\beta(y_0,0) g^{\alpha\beta}_E \Big|^2  \\
&\qquad + {\cal R}(f,F,b,B)^2 \Bigg ) t^2 \label{eq:itotaylor1}
\end{split}
\end{equation}
where ${\cal R}(f,F,b,B)$ is a term independent of $a$, $A$ and $t$.
\end{proposition}
\begin{proof}
As an example of how to compute the operators $L_\xi$
for the system of
equations \eqref{xyEquations}, 
we write down $L_{(\alpha)}$.
\[
L_{(\alpha)} f = 
b^i_\alpha \frac{\partial f}{\partial x^i}
+ B^i_\alpha \frac{\partial f}{\partial y^i}.
\]
Let us now the first few terms of the It\^o--Taylor expansion
for $z=f(x)-F(y)$.
\[
\begin{split}
L_{(0)}(z) &= f_*( a(x_t,t)) - F_*(A(y_t,t)) \\
&\quad + \frac{1}{2}(\D_{b_\alpha(x_t,t)}f_*)b_\beta(x_t,t) g^{\alpha \beta}_E - \frac{1}{2}(\D_{B_\alpha(y_t,t)}F_*)B_\beta(y_t,t) g^{\alpha \beta}_E. \\
L_{(\alpha)}(z) &= f_*( b_\alpha( x_t,t)) - F_*(B_\alpha(y_t,t)). \\
L_{(\alpha,\beta)}(z)&=L_{(\alpha)}L_{(\beta)}(z) \\
&= b_\alpha^i(x_t,t) \frac{\partial }{\partial x^i}f_*( b_\beta( x_t,t)) - 
B_\alpha^i(y_t,t) \frac{\partial }{\partial y^i}F_*(B_\beta(y_t,t)).
\end{split}
\]
We can now write down the order $1$ It\^o-Taylor expansion $z^1_t$.
It is
\[
\begin{split}
z^1_t &= \Big( f_*( a(x_0,0)) - F_*(A(y_0,0)) \\
&\qquad + \frac{1}{2}(\D_{b_\alpha(x_0,0)}f_*)b_\beta(x_0,0) g^{\alpha \beta}_E
- \frac{1}{2}(\D_{B_\alpha(y_0,0)}F_*)B_\beta(y_0,0) g^{\alpha \beta}_E \Big) I^{(0)}_{0,t}  \\
&\quad + \Big(f_*( b_\alpha( x_0,0)) - F_*(B_\alpha(y_0,0))\Big) I^{(\alpha)}_{0,t} \\
&\quad + \left(
b_\alpha^i(x_0,0) \frac{\partial }{\partial x^i}f_*( b_\beta( x_0,0)) - 
B_\alpha^i(y_0,0) \frac{\partial }{\partial y^i}F_*(B_\beta(y_0,0))
\right) I^{(\alpha,\beta)}_{0,t}
\end{split}
\]
We can now use \Cref{momentCalculations} to calculate $E(|z^1_t|^2)$.
This gives the desired result.
\end{proof}

\begin{remark}
For readers familiar with the traditional It\^o formula in Euclidean spaces, the term $(\D_{b_{\alpha(x_0,0)}}f_{i,*}) b_\beta(x_0,0) g^{\alpha\beta}_E$ for the $i$-th component of $f$ might be more familiar when written as  
\[  (\D_{b_{\alpha}}f_{i,*}) b_\beta g^{\alpha\beta}_E = {\mbox{Tr}}\left[ b^T (H f_i) b \right] \]
where Tr is the trace operator and $H$ is the Hessian operator. 
\end{remark}

\begin{theorem}[It\^o--Taylor series and It\^o-vector projection]
Given any time $t>0$,
if we wish to find the coefficients $A$ and $B$ at time $0$ for which
the solution to \Cref{equationOnManifold} is as close as possible to
the solution to \Cref{equationOnRr} in the sense that the 
mean square ($L^2$) norm of the order
$\frac{1}{2}$ It\^o--Taylor series for $X_t-\phi(Y_t)$
is minimized, we must take

\begin{equation*}
B_\alpha(Y_0,0) = (\psi_*)_{X_0} \Pi_{X_0} b_\alpha(X_0,0) 
\end{equation*}
where $\Pi_{X_0}$ is the projection map onto the tangent space of $M$ at $X_0$.
If we now suppose that $B$ is
chosen so that this minimum
is achieved at all points of $U$, a neighborhood of $X_0$ in $M$,
then the mean square $L^2$ norm of the order
$1$ It\^o--Taylor series is minimized by taking 
\[
A(Y_0,0) = (\psi_*)_{X_0} \Pi_{X_0} \left(a(X_0,0) - \frac{1}{2}(\D_{B_\alpha(Y_0,0)} \phi_*) B_\beta(Y_0,0) g^{\alpha\beta}_E \right).
\]
\label{theoremStrongVector}
\end{theorem}
\begin{proof}
We apply \Cref{generalProposition} taking $f$ equal to the identity,
$F$ equal to $\phi$,
$x_t=X_t$ and $y_t=Y_t$. To minimize the order $\frac{1}{2}$
It\^o--Taylor series for $X_t-\phi(Y_t)$ we must solve the problem:
\[
\hbox{Find $B_\alpha(Y_0,0)$ minimizing }\sum_{\alpha}|\phi_*(B_\alpha(Y_0,0)) - b_\alpha(X_0,0)|^2.
\]
The solution to this is given by $B_\alpha(Y_0,0)=\psi_*(V_\alpha)$ where 
the vectors $V_\alpha$ give a solution to the problem:
\[
\hbox{Find $V_\alpha \in \Image \phi_*$ minimizing }\sum_{\alpha}|V - b_\alpha(X_0,0)|^2.
\]
The standard properties of the projection map tell us that $V_\alpha = \Pi_{X_0}
b_\alpha(X_0,0)$. 

The same argument is used to find the formula for the coefficient $A$
that minimizes the order $1$ It\^o--Taylor expansion.
\end{proof}

This theorem motivates the following definition.

\begin{definition}
The It\^o-vector projection of the SDE \eqref{equationOnRr} onto
the manifold $M$ is given in the chart $\psi$ by the
SDE \eqref{equationOnManifold} with 
\begin{equation}
\begin{split}
\phi&:=\psi^{-1} \\
B_\alpha(Y_t,t) &:= (\psi_*)_{\phi (Y_t)} \Pi_{\phi (Y_t)} b_\alpha(\phi (Y_t),t) \\
A(Y_t,t) &:= (\psi_*)_{\phi (Y_t)} \Pi_{\phi (Y_t)} \left(a(\phi(Y_t),t) - \frac{1}{2}(\D_{B_\alpha(Y_t,t)} \phi_*) B_\beta(\phi(Y_t),t) g^{\alpha\beta}_E \right)
\end{split}
\label{itoProjection}
\end{equation}
\end{definition}

\begin{remark} The optimal $B$ in the above definition is the same we had in the Stratonovich projection in Eq. \eqref{stratProjectionFull}. The optimal $A$ is different. 
\end{remark}

\begin{corollary}
\label{independenceCorollary}
The It\^o-vector projection defines an SDE on the manifold $M$. By this we
mean that SDE defined on the manifold $M$ transforms
according to It\^o's lemma as we change chart $\psi$. See
\cite{armstrongBrigoJets} for a more detailed discussion of
the It\^o formulation of SDEs on manifolds.
\end{corollary}
\begin{proof}
The criteria we are using for finding the optimal coefficients of the SDE is given in terms of an estimate of the growth of the difference between the solution to the SDE in $\R^r$ and the solution to the
SDE on the manifold. Since it is expressed in terms of the solutions
to the SDE rather than the coefficients of the SDE, the criterion is independent of the choice of chart $\psi$. 

It follows that the condition we have derived on the coefficients will transform according to It\^o's lemma as we change the choice of chart. For an alternative proof by brute-force calculation see \cite{armstrongBrigoicms}.
\end{proof}

We will demonstrate that the It\^o-vector projection is
distinct from the Stratonovich projection by calculating an explicit
example in \Cref{applicationSection}.

One criticism of our derivation of the It\^o-vector projection
is that it is peculiar to worry about minimizing a term of
order $1$ when we cannot even ensure that the projection is accurate
to order $\frac{1}{2}$. It seems uncontroversial that
choosing the diffusion coefficient $B$ by the prescription above will yield
the best approximation, but will it make much difference to
choose $A$ in the same way?

The reason that choosing $A$ is
important
is that the errors of order $\frac{1}{2}$ due to the approximation of
$b$ will cancel on average. The correct choice of $A$ yields the optimal average value for the approximation. This is made precise by
the next result.

\begin{theorem}[It\^o-vector projection and weak It\^o-Taylor expansion]
If we wish to find the coefficient $A$ at time $0$ for which
the solution to \Cref{equationOnManifold} is as close as possible to
the solution to \Cref{equationOnRr} in the sense that the 
norm of the expectation of the order 
$1$ weak It\^o--Taylor series for $X_t-\phi(Y_t)$
is minimized, we must take:
\[
A(Y_0,0) = (\psi_*)_{X_0} \Pi_{X_0} \left(a(X_0,0) - \frac{1}{2}(\D_{B_\alpha(Y_0,0)} \phi_*) B_\beta(Y_0,0) g^{\alpha\beta}_E \right).
\]
where $\phi=\psi^{-1}$.
\label{theoremWeakVector}
\end{theorem}
\begin{proof}
The expectation of the weak It\^o Taylor expansion of $X_t-\phi(Y_t)$ is
\[
\left(a(X_0,0) - \phi_*(A(Y_0,0)) 
- \frac{1}{2}(\D_{B_\alpha(Y_0,0)} \phi_*) B_\beta(Y_0,0) g^{\alpha\beta}_E
\right) t.
\]
The result now follows immediately from the properties of $\Pi$.
\end{proof}

Thus the It\^o-vector projection is the choice
of $A$ and $B$ that simultaneously minimizes the expectation of the
error to order $\frac{1}{2}$ and the error of the expectation to order $1$.

Given the Stratonovich--Taylor expansions described in \cite{kloedenAndPlaten}, one might wonder if there are versions of 
Theorems \ref{theoremStrongVector} and \ref{theoremWeakVector} using Stratonovich--Taylor
series in place of It\^o--Taylor series? Might these provide a justification for the Stratonovich projection? The answer is negative.
The order 1
Stratonovich--Taylor expansion of \cite{kloedenAndPlaten} is in fact equal
to the order 1 It\^o--Taylor expansion. The difference is simply that
the Stratonovich--Taylor expansion is expressed in terms of Stratonovich coefficients and Stratonovich integrals rather than It\^o coefficients and integrals. Thus there is no different ``Stratonovich'' version of \Cref{theoremStrongVector}.

However, there is a sense in which the Stratonovich projection is optimal in relation with time symmetry. We will  address this optimality after introducing two different optimal approximations, the It\^o vector and It\^o jet projections.

\subsection{It\^o-jet projection}

\begin{figure}[h!]
\center{\includegraphics{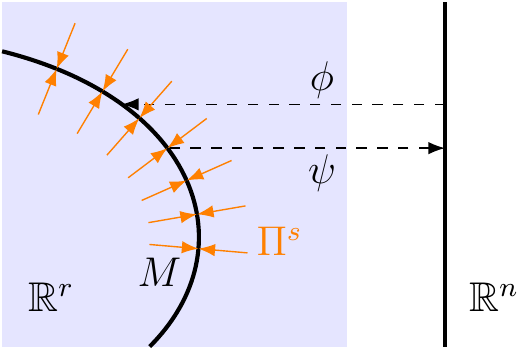}}
\caption{Metric projection $\Pi^s$ of a tubular neighbourhood of $M$ in $\R^r$
onto a neighbourhood $U$ in $M$. This is used to define the It\^o-jet projection.} 
\end{figure}

We now suppose that the open set $U$ inside our manifold $M$
has been chosen so that we can find a tubular neighbourhood $N$
of $U$ such that the metric projection $\Pi^s$ is smoothly
defined on $N$. The metric projection is the map sending a 
point $x \in \R^r$ to the nearest point in $N$. The standard
theory of tubular neighbourhoods tells us that if we choose $U$
small enough, these conditions will apply.

Note that the superscript $s$ in $\Pi^s$ is short for smooth and is intended
to distinguish this map from the linear projection operator $\Pi_{\phi(x)}$ onto
the tangent space at $\phi(x)$ ($x \in \R^n$).
Since the metric on $M$ is induced by the $\R^r$ Euclidean metric, we will have that the tangent-space linear projection, $\Pi_{\phi(x)}$, will be the first-order-component or best-linear-approximation of the metric projection, $\Pi^s$. See also our explicit calculation in \Cref{piFirstOrder} later on.

For ODEs only the first order linear component of the metric projection is
necessary to define the projection. However, It\^o SDEs involve explicit second order effects, so that there is an actual difference in applying the tangent vector projection or the full metric projection, going beyond the linear term, in approximating a SDE on a submanifold. As we pointed out in \cite{armstrongBrigoJets}, an It\^o SDE can be interpreted as a 2-jet. It is then not completely surprising that the second order terms of the metric projection play an important role in understanding the projection of SDEs.

More specifically, in this section we will solve the problem of finding an
SDE  on the manifold $M$, $Y_t$ in $\psi$ coordinates, which
minimizes the mean square of the truncated Taylor expansion of the $M$ geodesic distance between $\Pi^s(X_t)$ and
$\phi(Y_t)$, or ambient $\R^r$ distance between these two points of $M$. The two distances will lead to the same result. We call this solution the It\^o-jet projection.
By contrast, the It\^o-vector projection focuses on the $\R^r$ distance between $\phi(Y_t)$ and $X_t$.
Thus the It\^o-jet projection uses the metric projection of $X$ as a benchmark to obtain an optimal approximation $\phi(Y)$, whereas the It\^o-vector projection uses directly the original $X$ as a benchmark.

The It\^o-jet projection is most neatly defined using the correspondence
between $2$-jets and SDEs described in \cite{armstrongBrigoJets}. We
recall this correspondence now.

Suppose that at each point, $x$, of a manifold, $M$, we are
given a smooth map $\gamma_x:\R^m \to M$
with $\gamma_x(0)=x$. Suppose also that $\gamma_x$
depends smoothly on $x$. We can then define an SDE on $M$ driven
by $m$ dimensional Brownian motion by:
\[
\ed X^i_t = \frac{1}{2} g^{\alpha \beta}_E \frac{\partial^2 \gamma^i}{\partial x^\alpha \partial x^\beta}(X_t) dt
+ \frac{\partial \gamma^i}{\partial x^\alpha}(X_t) dW^\alpha_t .
\]
Here $X^i$ and $\gamma^i$ denote the components in some coordinate chart on $M$. It follows by It\^o's Lemma that this SDE is independent of the
choice of charts for $M$. Since we only use the first two derivatives
of $\gamma$ in this definition, we say that the SDE depends only on the
$2$-jet of $\gamma$.

We
can now write down the definition of the It\^o-jet projection.

\begin{definition}
Let $W^\alpha_t$ be independent Brownian motions with $1 \leq \alpha \leq k$.
Let $\gamma_x:\R^k \to \R^r$ be a smoothly varying family of maps satisfying $\gamma_x(0)=x$ for all $x \in \R^r$. We interpret $\gamma$ as defining an It\^o SDE. We define the It\^o-jet projection
to be the SDE associated with $\Pi^s \circ \gamma_y:\R^k \to M$. 
\end{definition}

Since this
definition only depends upon germs of $\Pi^s$ and $\gamma$, the It\^o-jet projection does not depend upon issues such as the tubular neighbourhood used 
to define $\Pi^s$.

We wish to show that the It\^o-jet projection solves
the problem of finding the best approximation to the SDE on the manifold,
if one measures the quality of the approximation using the truncated It\^o-Taylor expansion of either the 
 geodesic distance or the  distance in the ambient space $\R^r$.

\begin{theorem}[It\^o-jet projection as optimal approximation]
\label{jetTheorem}
Let $\lambda(x,y)$ denote the square of the geodesic distance between two points
on $M$. Let $|x-y|_r^2$ denote the square of the distance between two points in
the ambient space.
If we wish to find the coefficients $A$ and $B$
at time 0 for which the solution to \eqref{equationOnManifold}
is as close as possible to the image on $M$ of the solution of \eqref{equationOnRr} 
under $\Pi^s$ in the sense that
the expectation of the square of the order $1$ 
It\^o-Taylor expansion of $\lambda$ or of $|\cdot|_r^2$ is minimized, we
must take
\begin{equation*}
B_\alpha(Y_0,0) = (\psi_*)_{X_0} \Pi_{X_0} b_\alpha(X_0,0).
\end{equation*}
If we use this to define $B$ at all points of $M$, we have
that the expectation of the square of the order $2$
It\^o-Taylor expansion of  $\lambda$ or $|\cdot|_r^2$ is minimized 
by ensuring that the the $2$-jet associated with 
\eqref{equationOnManifold} at $(Y_0,0)$ is given
by $\Pi^s \circ \gamma_{X_0}$ where $\gamma_x$ is the
$2$-jet associated with \eqref{equationOnRr}. This results in the following drift
for \eqref{equationOnManifold}:
\begin{equation}\label{eq:Ajettheorem}
A(Y_0,0)= \tilde{\Pi}^s_*(a(X_0,0))
+ \frac{1}{2} (\D_{b_{\alpha(X_0,0)}}\tilde{\Pi}^s_*)b_\beta(X_0,0) g^{\alpha\beta}_E,
\end{equation}
where we define $\tilde{\Pi}^s=\psi \circ \Pi^s$.
\end{theorem}
\begin{proof}
We will first prove the result for the geodesic distance.

It will suffice to prove the result in a single chart. Hence
we may assume that our coordinates are normal coordinates based
at $X_0$.

We have the following Taylor
series expansion for the square of the geodesic distance (see for example formula 3.4.3 in \cite{brewin}):
\[
\begin{split}
\lambda(x,y) &= g^E_{ij}(x^i - y^i)(x^j - y^j) \\
&\quad - \frac{1}{12}  R_{ikjl}(x^i+y^i)(x^k-y^k)(x^j+y^j)(x^l-y^l) + O((|x|+|y|)^5).
\end{split}
\]
The first term is just the Euclidean metric on $\R^n$, the term
$R_{ikjl}$ denotes the Riemann curvature tensor of $M$ at the origin.

We can
write down the expectation of the order $2$
It\^o--Taylor expansion of $|\tilde{\Pi}^s(X_t)-Y_t|^2$
using \Cref{generalProposition}, taking $f=\tilde{\Pi}^s$ and $F$ to be the identity. It is:
\begin{equation}
\begin{split}
&\sum_{\alpha}|\tilde{\Pi}^s_*(b_\alpha(X_0,0))- B_\alpha(Y_0,0)|^2 t \\
&\quad + \Bigg ( \left|\tilde{\Pi}^s_*(a(X_0,0))
+ \frac{1}{2} (\D_{b_{\alpha(X_0,0)}}\tilde{\Pi}^s_*)b_\beta(X_0,0) g^{\alpha\beta}_E
- A(Y_0,0) \right|^2  \\
&\qquad + {\cal R}(\tilde{\Pi}^s,b,B)^2 \Bigg ) t^2
\end{split}
\label{order2ItoTaylorSeries}
\end{equation}
where ${\cal R}(\tilde{\Pi}^s,b,B)$ is a term independent of $a$ and $A$. 
Our reasoning is that
we know that this formula for the expectation of $|\tilde{\Pi}^s(X_t)-Y_t|^2$
is accurate up to order $t^2$, therefore it must equal the 
expectation of the order $2$ It\^o--Taylor expansion. This allows us to avoid
computing the order $2$ It\^o--Taylor expansion directly.

The curvature term is fourth order, so it will not influence the
order $1$ It\^o--Taylor expansion for $\lambda$. This is because the differential operators $L_\xi$ in this expansion are all order $2$ or less. We deduce that
the expectation of the order 1 It\^o-Taylor expansion of $\lambda$
is
\[
\sum_{\alpha}|\tilde{\Pi}^s_*(b_\alpha(X_0,0))- B_\alpha(Y_0,0)|^2 t.
\]
This is minimized by taking $B$ as described in the statement of the theorem.

The expectation of an integral $I_\xi(1)$ is zero if $\xi$
contains any non-zero entries. This follows by the Martingale property
of the It\^o integral. Thus the non-zero terms in the expectation
of the order $2$ It\^o expansion for $\lambda$ correspond to
the multi-indices $()$, $(0)$ and $(0,0)$. Since the curvature term
is fourth order, the only term that will contain a curvature term
corresponds to the index $(0,0)$. Moreover, only the highest
order term of the operator $L_{0,0}$ is influenced by the curvature.
The coefficient of this highest order term may involve only $b$ and $B$
but will not involve $a$ or $A$.

Thus the expectation of the order $2$ It\^o--Taylor
expansion is of the form \eqref{order2ItoTaylorSeries} since
any curvature correction can be absorbed into the term 
${\cal R}(\tilde{\Pi}^s,b,B)^2$. We deduce that the order $2$ It\^o--Taylor
series is minimized by taking $A$ as in Equation \eqref{eq:Ajettheorem}.
When these conditions are rewritten in the language of $2$-jets, we
get the desired result for the metric $\lambda$.

The proof for the metric $|\cdot|$ follows from \Cref{lemma:ambientMetric} given below, and is otherwise essentially identical to that for $\lambda$.
\end{proof}

Note that in this argument we can ensure that the order $1$ expansion of
$\lambda$ actually vanishes. By contrast, recall that
the corresponding term did not vanish in the derivation of the
It\^o-vector projection which lead us to give an alternative
derivation using the weak It\^o--Taylor expansion.

\begin{lemma}
\label{lemma:ambientMetric}
Let $U$ be a neighbourhood of the origin in $\R^n$ and let $\phi:U \to \R^r$ be normal coordinates for the Riemannian manifold $\phi(U)$ centred at the origin, then
\[
|\phi(x)-\phi(y)|_r^2 = |x-y|^2_n + O((|x|_n+|y|_n)^4).
\]
Here $|\cdot|_n$ is the norm on $\R^n$.
\end{lemma}
\begin{proof}
Without loss of generality we may assume that the origin is mapped to the origin and the coordinate axes in $\R^n$ are
mapped to the corresponding axes in $\R^r$. Given a point $y\in U$, we can write the Taylor expansion for the component $\phi(y)^a$ in the following form:
\begin{equation}
\phi(y)^a = (\delta^r_n)^a_i y^i + A\indices{^a_{jk}} y^j y^k + O(|y|_n^3).
\label{eq:phiGeodesicNormalCoordinates}
\end{equation}
Here $(\delta^r_n)$ is the tensor representing the projection of $\R^r$ onto $\R^n$. The upper indices of $(\delta^r_n)$ range from 1 to $r$ and the lower from $1$ to $n$. $(\delta^r_n)^i_j$ is equal to 1 if $i=j$ and $0$ otherwise. 
$A\indices{^a_{jk}}$ is a tensor with upper index $a$ ranging from $1$ to $r$ and lower indices $j$ and $k$ ranging from $1$ to $n$ and which satisfies $A\indices{^a_{jk}}=A\indices{^a_{kj}}$.

The components of the metric tensor on $U$ can now be computed as follows:
\[
\begin{split}
g_{ij} &= \left\langle \frac{\partial \phi}{\partial y^i}, \frac{\partial \phi}{\partial y^j} \right\rangle_r \\
&= \left((\delta^n_r)^a_i + A\indices{^a_{ik}}y^k\right)
   \left((\delta^n_r)^b_j + A\indices{^b_{jl}}y^l\right) (g^r)_{ab}
   + O(|y|_n^2)
\end{split}
\]
Here $g^r$ is the metric tensor of $\R^r$. Our expression for $g_{ij}$ simplifies to give:
\[
g_{ij} = (g^n)_{ij} + A\indices{^j_{ik}}y^k + A\indices{^i_{jk}}y^k  + O(|y|_n^2).
\]

It is well known that in Riemannian normal coordinates the partial derivatives
of the metric tensor vanish at the origin. We compute that
\[
\partial_k g_{ij} \rvert_0 = A\indices{^j_{ik}} + A\indices{^i_{jk}}.
\]
So we have $A\indices{^j_{ik}} = -A\indices{^i_{jk}}$. However, recall that $A\indices{^i_{jk}}$ is symmetric in the indices $j$ and $k$. We see that:
\[
A\indices{^i_{jk}} = -A\indices{^j_{ik}} = -A\indices{^j_{ki}} = A\indices{^k_{ji}} = A\indices{^k_{ij}} = -A\indices{^i_{kj}} = -A\indices{^i_{jk}}.
\]
So all the components of $A$ vanish.

We can now use \eqref{eq:phiGeodesicNormalCoordinates} to compute:
\[
\begin{split}
|\phi(x)-\phi(y)|_r^2 &= |(\delta^r_n)^a_i x^i - (\delta^r_n)^a_i y^i|_r^2 + O((|x|_n + |y|_n)^4) \\
&= |x-y|_n^2  + O((|x|_n + |y|_n)^4).
\end{split}
\]
\end{proof}

\subsection{Time-symmetric optimality of the Stratonovich projection}

Having introduced the It\^o vector and It\^o jet projections, we are now in a position where we can clarify that also the Stratonovich projection is optimal in a time symmetric sense, even if this optimality is somewhat ad hoc. 

Consider the $x$ SDE \eqref{xyEquations}, in Stratonovich form:
\begin{equation}
\ed x_t = \bar{a}(x,t) \ed t + b_\alpha(x,t) \circ \ed {W}^\alpha_t, \ \ X_0.
\label{generalSDEstrat}
\end{equation}
Recall the It\^o Stratonovich transformation 
\[ \bar{a}_i = a_i - \frac{1}{2} \sum_{j=1}^{d_x}  \sum_{\alpha} b_\alpha^j \frac{\partial b_\alpha^i}{\partial x_j} . \]

More generally, by a bar over the drift of an It\^o SDE we will mean the drift of the equivalent Stratonovich SDE. 

We now extend the SDE to negative time as follows. Define
\begin{equation}
\ed \xi_t = -\bar{a}(\xi,t) \ed t - b_\alpha(\xi,t) \circ \ed \hat{W}^\alpha_t, \ \ \xi_0 = X_0
\label{generalSDEstrat}
\end{equation}
where $\hat{W}$ is a second standard Brownian motion, independent of ${W}$.
Given the symmetric nature of the Stratonovich integral underlying the above SDE and given that formally the chain rule holds, it makes sense to define $x$ for $t<0$ by setting
\[ x_{-t} := \xi_t .\]


We now wonder whether the Stratonovich projection could be indeed optimal at time $0$ for this SDE extended to negative time at time $0$. Suppose that we wish to find the SDE on $M$
\begin{equation}
\ed y = \bar{A} \, \ed t + B_\alpha \circ \ed W^\alpha_t, \qquad y_0=x_0,
\label{equationOnManifoldstratop}
\end{equation}
extended similarly to negative time (giving $y_{-t}$), that minimizes the mean square of the truncated Taylor expansion of  the vector $(f(x_t)-F(y_t),f(x_{-t})-F(y_{-t}))$ 
Here $f$ and $F$ are functions as defined in Proposition \ref{generalProposition}.

This optimality criterion is symmetric under time reversal around an anchor state given by a deterministic initial condition. For most applications, for example, when we apply the projection method to filtering, there will be a clear time asymmetry in the problem setting. In these cases, a time-symmetric optimality criterion would not be appropriate. However, in applications to physics one may possibly seek to approximate SDEs that are symmetric under time-reversal in a manner that preserves this symmetry. In this case this criterion would be a natural choice.  

Write $z^i_t$ for the terms in the It\^o Taylor expansion of $z_t=f(x_t)-F(y_t)$. Proposition \ref{generalProposition} states that for positive $t$:
\begin{equation}
\begin{split}
E(|z^{\frac{1}{2}}_t|^2)
&= \sum_{\alpha}|f_*(b_\alpha(x_0,0))- F_*(B_\alpha(y_0,0))|^2 t \\
E(|z^1_t|^2)
&= \sum_{\alpha}|f_*(b_\alpha(x_0,0))- F_*(B_\alpha(y_0,0))|^2 t \\
&\quad + \Bigg ( \Big|f_*(a(x_0,0)) - F_*(A(y_0,0)) \\
&\quad + \frac{1}{2} (\D_{b_{\alpha(x_0,0)}}f_*)b_\beta(x_0,0) g^{\alpha\beta}_E
 - \frac{1}{2} (\D_{B_{\alpha(y_0,0)}}F_*)B_\beta(y_0,0) g^{\alpha\beta}_E \Big|^2  \\
&\qquad + {\cal R}(f,F,b,B)^2 \Bigg ) t^2 \label{eq:itotaylor1}
\end{split}
\end{equation}
The interesting point is that for negative $t$ the correction terms in
this second equation are very different. Recall $\xi_t=x_{-t}$ and define $\eta_t=y_{-t}$ so we can write conventional
It\^o SDEs for $\xi_t$ and $\eta_t$.
\[
\ed \xi_t = -(a(\xi_t,-t)-(\D_{b_\alpha(\xi_t,-t)} b_\beta)(\xi_t,-t) g^{\alpha \beta}_E) \ed t - b_\alpha(\xi_t,-t) \ed \hat{W}^\alpha_t 
\]
\[
\ed \eta_t = -(A(\eta_t,-t)-(\D_{B_\alpha(\eta_t,-t)} B_\beta)(\eta_t,-t) g^{\alpha \beta}_E) \ed t - B_\alpha(\eta_t,-t) \ed \hat{W}^\alpha_t  .
\]
By proposition \ref{generalProposition} we therefore have the following expressions for $z^i_t$
for negative times $t$.
\begin{equation}
\begin{split}
E(|z^{\frac{1}{2}}_t|^2)
&= \sum_{\alpha}|f_*(b_\alpha(x_0,0))- F_*(B_\alpha(y_0,0))|^2 t \\
E(|z^1_t|^2)
&= \sum_{\alpha}|f_*(b_\alpha(x_0,0))- F_*(B_\alpha(y_0,0))|^2 t \\
&\quad + \Bigg ( \Big|f_*\left(-a(x_0,0) + \D_{b_\alpha(x_0,0)}b_\beta(x_0,0) g^{\alpha \beta}_E \right) \\
&\quad - F_*\left(-A(y_0,0)  + \D_{b_\alpha(x_0,0)}b_\beta(x_0,0) g^{\alpha \beta}_E\right) \\
&\quad + \frac{1}{2} (\D_{b_{\alpha(x_0,0)}}f_*)b_\beta(x_0,0) g^{\alpha\beta}_E
 - \frac{1}{2} (\D_{B_{\alpha(y_0,0)}}F_*)B_\beta(y_0,0) g^{\alpha\beta}_E \Big|^2  \\
&\qquad + {\cal R^\prime}(f,F,b,B)^2 \Bigg ) t^2
\end{split}
\end{equation}

Combining our results we have for any $t$:
\[
\frac{1}{2} E(|z_t^{\frac{1}{2}}|^2
+ \frac{1}{2} E(|z_{-t}^{\frac{1}{2}}|^2) = \sum_{\alpha}|f_*(b_\alpha(x_0,0))- F_*(B_\alpha(y_0,0))|^2 t
\]
\begin{equation}
\begin{split}
\frac{1}{2} E(|z^1_t|^2) + \frac{1}{2} E(|z^1_{-t}|^2)
&= \sum_{\alpha}|f_*(b_\alpha(x_0,0))- F_*(B_\alpha(y_0,0))|^2 t \\
&\quad + \Bigg ( \frac{1}{2} \Big|f_*\left(-a(x_0,0) + \D_{b_\alpha(x_0,0)}b_\beta(x_0,0) g^{\alpha \beta}_E \right) \\
&\quad \quad - F_*\left(-A(y_0,0)  +  \D_{B_\alpha(x_0,0)}B_\beta(x_0,0) g^{\alpha \beta}_E\right) \\
&\quad \quad + \frac{1}{2} (\D_{b_{\alpha(x_0,0)}}f_*)b_\beta(x_0,0) g^{\alpha\beta}_E
 - \frac{1}{2} (\D_{B_{\alpha(y_0,0)}}F_*)B_\beta(y_0,0) g^{\alpha\beta}_E \Big|^2  \\
&\quad + \frac{1}{2} \Big|f_* a(x_0,0) - F_* A(y_0,0) \\
&\quad \quad + \frac{1}{2} (\D_{b_{\alpha(x_0,0)}}f_*)b_\beta(x_0,0) g^{\alpha\beta}_E
  - \frac{1}{2} (\D_{B_{\alpha(y_0,0)}}F_*)B_\beta(y_0,0) g^{\alpha\beta}_E \Big|^2  \\ 
&\quad + {\cal R}(f,F,b,B)^2 + {\cal R^\prime}(f,F,b,B)^2 \Bigg ) t^2
\label{eq:itoTaylor1Symmetric}
\end{split}
\end{equation}

We can now find the optimal projection in a time symmetric sense by mimicking our previous arguments. We will need a lemma that allows us to write down the minimizing $A$.
\begin{lemma}
\label{lemma:optimization}
Let $(V,g)$ be a Hilbert space containing a closed subspace $W$.
Let $v_1$ and $v_2$ be two vectors in $V$. Then the optimization problem
\[
\begin{aligned}
\underset{w\in W }{\mathrm{minimize}} & |w-v_1|^2 + |w-v_2|^2
\end{aligned}
\]
has a unique minimizer given by $w = \Pi_W(\frac{1}{2} v_1 + \frac{1}{2}v_2)$
where $\Pi_W$ is orthogonal projection onto $W$.
\end{lemma}
\begin{proof}
By translating the problem by $\Pi_W(\frac{1}{2} v_1 + \frac{1}{2}v_2)$ we may assume wlog that $\Pi_W(v_1)=-\Pi_W(v_2)$. Then for $w \in W$
\[
\begin{split}
|w-v_1|^2 + |w-v_2|^2 &= |w-\Pi_W v_1|^2 + |w- \Pi_W v_2|^2 + |\Pi_W v_1 - v_1|^2
 + |\Pi_W v_2 - v_2|^2 \\
 &= |w-\Pi_W v_1|^2 + |w+\Pi_W v_1|^2 + |\Pi_W v_1 - v_1|^2
  + |\Pi_W v_2 - v_2|^2 \\
 &= 2|w|^2 + 2|\Pi_W v_1|^2 + |\Pi_W v_1 - v_1|^2
  + |\Pi_W v_2 - v_2|^2  
 \end{split}
\]
This has a unique minimizer, $w=0$.
\end{proof}

Let us begin with the case $f=\mathrm{id}$ and $F=\phi$. Minimizing the order $\frac{1}{2}$ expansion requires us to choose $B$ such that
$\Pi_{\phi(y_t)} b_\alpha(\phi(y_t),t)=\phi_*(B_\alpha(y_t),t)$ at every point.
By Lemma \ref{lemma:optimization} and \eqref{eq:itoTaylor1Symmetric}, minimizing the order $1$ expansion then requires us to choose $A$ such
that
\[
\phi_*(A(y_t,t))=\Pi_{\phi(y_t)} \left( a(\phi(y),t)
- \frac{1}{2} \D_{b_\alpha(\phi(y_t),t)} b_\beta(\phi(y_t),t) g^{\alpha \beta}_E
+ \frac{1}{2} \phi_* \D_{B_\alpha(\phi(y_t),t)} B_\beta(\phi(y_t),t) g^{\alpha \beta}_E \right).
\]
Equivalently
\[
\phi_*(\bar{A}(y_t,t))=\Pi_{\phi(y_t)} \left( \bar{a}(y,t) \right).
\]
Hence we obtain the Stratonovich projection as the time symmetric analogue of the vector projection.

Now consider the case where $f=\Pi_s$ and $F=\mathrm{id}$. Minimizing the order 
$\frac{1}{2}$ expansion requires us to choose $B$ by
$(\Pi_s)_* b_\alpha = B_\alpha$. By Lemma \ref{lemma:optimization}
and \eqref{eq:itoTaylor1Symmetric}, minimizing the order $1$ expansion then requires us to choose $A$ such
that
\[
A(y_t,t)=(\Pi_s)* \left( a(\phi(y),t)
- \frac{1}{2} \D_{b_\alpha(\phi(y_t),t)} b_\beta(\phi(y_t),t) g^{\alpha \beta}_E \right)
+ \frac{1}{2} \D_{B_\alpha(\phi(y_t),t)} B_\beta(\phi(y_t),t) g^{\alpha \beta}_E .
\]
Equivalently
\[
\bar{A}(y_t,t)=(\Pi_s)_* \left( \bar{a}(\phi(y),t) \right).
\]
Again, this is the Stratonovich projection.

\section{A low dimensional example: cross diffusion on a unit circle}
\label{section:lowDimensionalExample}

We now look at a concrete example which shows the difference between the It\^o-vector and It\^o-jet projections. Consider the SDE in $\R^2$ given by
\begin{equation}
\label{crossDiffusion}
\begin{split}
\ed X_t &= \sigma Y_t \ dW_t, \\
\ed Y_t &= \sigma X_t \ dW_t,
\end{split}
\end{equation}
with deterministic initial condition $(X_0,Y_0)$. We call this a cross diffusion, since each state crosses over as diffusion coefficient of the other state and the paths tend to lie on a St Andrew cross, see \Cref{appendixcross} for more details on this process.
We wish to project this process equation onto the unit circle given by $X^2+Y^2=1$.
It is easy to check using It\^o's Lemma that if we write $(X_t,Y_t)$ in polar
coordinates as $(r_t \cos(\theta_t), r_t \sin(\theta_t))$ then $\theta_t = \arctan(Y_t/X_t)$
satisfies the following exact angular position process equation:
\begin{equation}
\ed \theta_t = -\frac{1}{2} \sigma^2 \sin(4 \theta_t) \ed t
		     + \sigma \cos(2 \theta_t) \ed W_t, \ \  \mbox{or} \ \ \  \ed \theta_t =  \sigma \cos(2 \theta_t) \circ \ed W_t.
\label{examplePolarEquation}		     
\end{equation}

Thanks to the special structure of the cross-diffusion, the equation above is already a closed SDE for $\theta$ without needing to apply any of our projection methods. In this sense we already have the exact angular position SDE and we do not need to project the original $\R^2$ SDE on the circle $M$ to approximate the exact angular position with a SDE on the circle. However, we might want to check whether one of our projection methods is consistent with the exact angular position SDE. Let us check how the different projections behave. 
If we use the same polar coordinate $\theta$ for the unit circle, we find
that the Stratonovich projection and the It\^o-jet projection for the $(X,Y)$ SDE are also given by \eqref{examplePolarEquation}, and are thus consistent with the exact $\theta$. However the It\^o-vector
projection is different and results in:
\begin{equation*}
\ed \theta_t = \sigma \cos(2 \theta_t) \ed W_t.
\end{equation*}

For this example at least, the It\^o-jet projection and the
Stratonovich projections track the angular position of $(X_t,Y_t)$
perfectly. Intuitively one might therefore feel that the Stratonovich
and  It\^o-jet projections are ``better'' approximations to the SDE
despite the short time optimality arguments given earlier. It turns out this is a special case of a more general situation, summarized in the following
\begin{definition}[SDE that fibers over a map between manifolds]
Let $f:M\to N$ be a smooth map between two manifolds. Let $S$
be an SDE on $M$ determined by the 2-jets $\gamma_x:\R^m\to M$
given at each point $x \in M$. We say that {\em $S$ fibres over $f$} if $j_2(f \circ \gamma_{x_1})=j_2(f \circ \gamma_{x_2})$
whenever
$f(x_1)=f(x_2)$. This implies that we can define an SDE on the image of $f$
using the $2$-jets $j_2(f \circ \gamma_x)$ at $f(x)$. We call this the
{\em SDE induced by $f$}.
\end{definition}

Returning to projection, we see that we have the following

\begin{theorem}[If SDE fibres over $\Pi^s$  then Stratonovich $=$ It\^o-jet proj.] If an SDE fibres over the smooth projection map $\Pi^s$ then the Stratonovich and It\^o-jet projection
will both be equal to the SDE induced by $\Pi^s$.
\end{theorem}
\begin{proof}
 This is an immediate consequence of the Stratonovich chain rule in the first case. It is a trivial
consequence of the definition of the It\^o-jet projection in the second case. 
\end{proof}
Our two-dimensional example of the cross-diffusion on the circle is simply a special case of this more general
phenomenon.

\begin{figure}[tbh]
\centering
\includegraphics[width=0.7\linewidth]{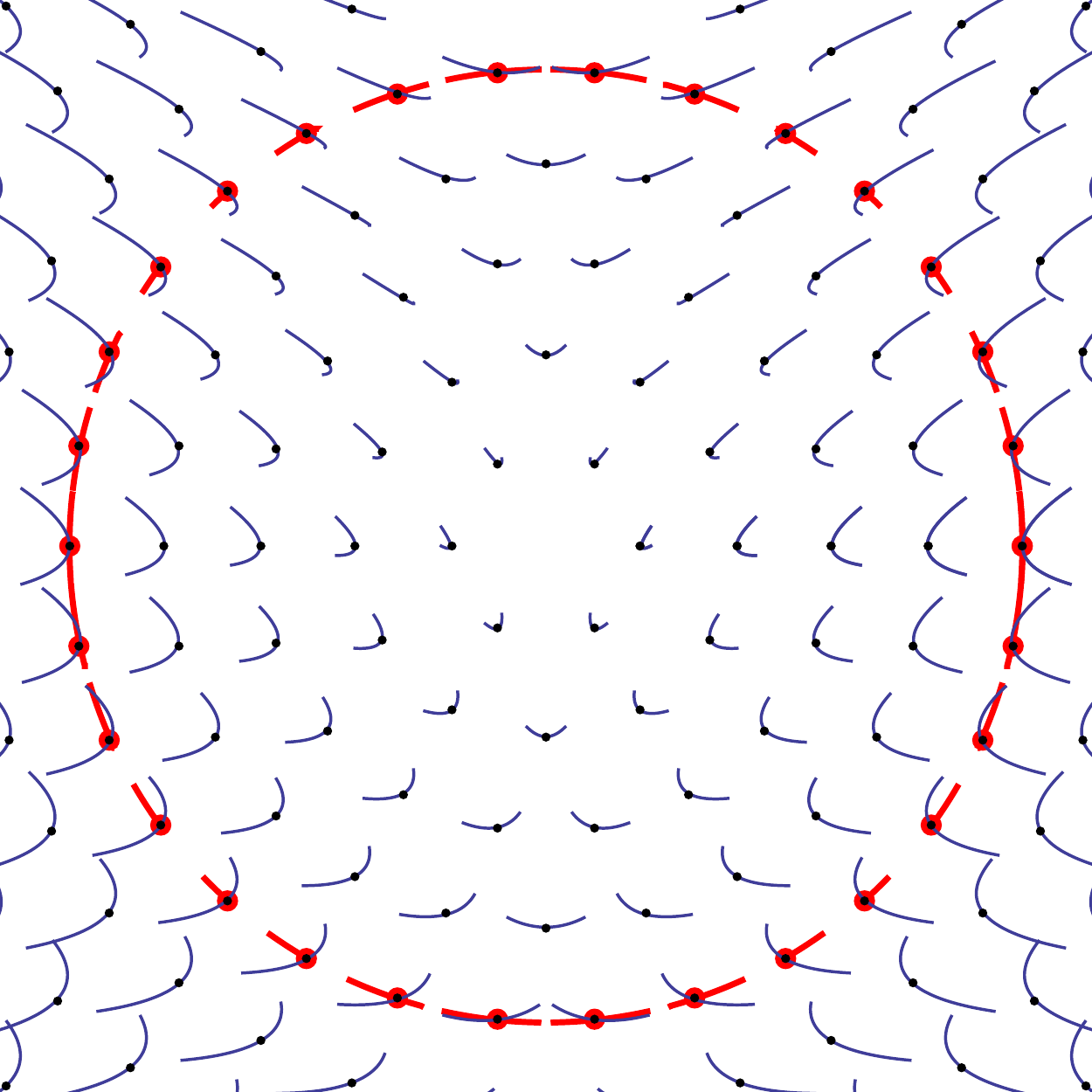}
\caption{An SDE in $\R^2$ and its It\^o-jet projection onto
the unit circle}
\label{fig:projectedSDE}
\end{figure}

It is interesting to note that one can draw a diagram to show
the It\^o-jet projection. In \cite{armstrongBrigoJets} it is discussed how
the jet formulation of SDEs makes it possible to draw pictures of SDEs
that transform according to It\^o's lemma. For processes driven by
one dimensional Brownian motion, one simply finds functions $\gamma_x$
whose $2$-jet represents the SDE and then draws the image of an interval
$[-\epsilon, \epsilon]$ under the map $\gamma_x$ at each point $x$.
A picture of this type is shown in \Cref{fig:projectedSDE}. It shows
how the $2$-jets determining the SDE \eqref{crossDiffusion} can be projected
onto the unit circle simply by composition with $\Pi^s$.

It seems paradoxical that we derived the It\^o-vector projection using optimality arguments that seem to be less ad hoc than for the Stratonovich projection, and yet, for this example, the
It\^o-vector projection appears manifestly suboptimal.

One possible resolution to this paradox is to say that our notions of tracking
$X_t$ optimally are flawed. \Cref{theoremStrongVector} has
the weakness that we attempt to minimize a term of order $1$ when our approximation is not accurate at order $\frac{1}{2}$.   Indeed, looking at \Cref{eq:itotaylor1} we see that when we try to minimize the relevant expectation we minimize a combination of terms of order $t$ and $t^2$ for the square. Moreover,  \Cref{theoremWeakVector} has the weakness that we are using the error in the mean to measure the accuracy of our solution. By contrast, the It\^o-jet projection has a fully convincing derivation as the optimal approximation of $\Pi^s(X_t)$ up to order 1.

We will see numerical evidence later that suggests that the 
It\^o-jet projection performs better in the long term than the It\^o-vector projection
which lends some support to the idea that the It\^o-jet projection is the ``right'' choice.

We summarize the different projections and the optimality criteria used to determine their drifts in \Cref{tab:projtypes}. The diffusion coefficient is identical for all three projections.

\begin{table}[tbhp]
\begin{center}
\begin{tabular}{p{20mm}p{80mm}}\toprule
Projection  & Properties of drift term  \\ \midrule
It\^o-vector
&

Minimizes norm of
the expectation of the order 1  weak It\^o--Taylor expansion between $X$ \& $\phi(Y)$.

(ii)
Given $B$ minimizing mean square of the order 1/2 
 strong It\^o--Taylor expansion 
of the difference $X_t-\phi(Y_t)$,
minimizes mean square of order 1 strong It\^o--Taylor expansion 
of the difference.
\\ \addlinespace                     
 It\^o-jet
  &
Minimizes mean square of $2$ Strong It\^o--Taylor 
 expansion for $\R^r$ or $M$   
 distance between $\Pi_s(X)$ \& $\phi(Y)$

\\ \addlinespace                     

Stratonovich
&
Similar to It\^o vector below but for the Taylor series of the differences vector $[X_t-\phi(Y_t), \ X_{-t}-\phi(Y_{-t})] $ at positive and negative time, where negative time processes are defined ad hoc by propagating a second input Brownian motion backward in time. 
 \\
                      \bottomrule
\end{tabular}
\end{center}
\caption{ Projections and the associated optimality criteria}\label{tab:projtypes}
\end{table}

\section{The It\^o-jet projection in local coordinates}
\label{section:secondOrderProjection}

Our definition of the It\^o-jet projection is coordinate free and
 simple. However, to calculate
it in practice we will need an explicit coordinate representation.

We therefore wish to calculate the metric projection map ${\tilde \Pi}^s=\psi \circ \Pi^s$ up to second order. Then using It\^o's formula for $2$-jets
we will be able to calculate the It\^o-jet projection associated to ${\tilde \Pi}^s$.

Most of our calculation involves the deterministic map ${\tilde \Pi}^s$. Thus
in this section we will drop the convention of using Greek indices
exclusively for components of the Brownian motion. In this section we
will also use Greek indices to highlight indices over which we are summing.
This makes the formulae a little easier to read.

We define the metric tensor on $U$ by:
\begin{equation}
h_{a b} = \frac{\partial \phi^\alpha}{\partial x^a}
    \frac{\partial \phi^\alpha}{\partial x^b}
    \label{metricTensor}
\end{equation}    
The differential ${\tilde \Pi}^s_*$ of ${\tilde \Pi}^s$ is well known to be given
by the linear projection onto $\Image \phi*$ composed with the map $\phi_*^{-1}$. Hence
${\tilde \Pi}^s_*$ is the unique linear map with ${\tilde \Pi}^s_* \circ \phi_*$ equal to the identity and with
kernel equal to the orthogonal complement of $\Image \phi*$. We deduce that ${\tilde \Pi}^s_*$ has
the following components:

\begin{equation}
 \Pi^a_b:= ({\tilde \Pi}^s_*)^a_b =  \frac{\partial \phi^b}{\partial x^\alpha} h^{a \alpha}, \qquad a \leq n, \alpha \leq n, b \leq r.
\label{piFirstOrder}
\end{equation}

We note that the differential or tangent map ${\tilde \Pi}^s_*$ is the best linear approximation of the metric projection ${\tilde \Pi}^s$ around the relevant point $x = \phi(y) \in M$, and it coincides with the classic linear projection $\Pi_{\phi(y)}$ on the tangent space of $M$. Indeed, \Cref{piFirstOrder} shows the classic components of the projection on the tangent space of an $n$-dimensional manifold $M$ embedded in $\R^r$ and realized as $\phi$-image of a subset or $\R^n$.

\begin{lemma}
\label{secondOrderProjectionLemma}
Suppose for simplicity that $\phi(0)=0$ and
\begin{equation*}
(\phi_*)^a_b := \frac{\partial \phi^a}{\partial x^b} = D^a_b :=
\begin{cases}
1 & a = b \hbox{ and } a \leq n\\
0 & \hbox{otherwise}
\end{cases}
\end{equation*}
then $\Pi^s$ is given up to second order by
\begin{equation*}
{\tilde \Pi}^s(y)^a = y^a
- \frac{1}{2} \frac{ \partial^2 \phi^a}{ \partial x^\alpha x^\beta} y^\alpha y^\beta
+ \frac{ \partial^2 \phi^\gamma}{ \partial x^a x^\beta} (\perp)_{\gamma \alpha} y^\alpha y^\beta
+ O(|y|^3)
\end{equation*}
where we define
\begin{equation*}
(\perp)_{a b} = 
\begin{cases}
1 & a = b \hbox{ and } a > n \\
0 & \hbox{otherwise}.
\end{cases}
\end{equation*}
Note that we are using an extension of the Einstein summation convention to cover tensors where
some indices range from $1$ to $n$ and some from $1$ to $r$. Where an index appears twice, we
sum over the smaller range. Note also that we are working in a restricted
set of coordinate systems, so it no longer holds that all summed pairs of indices will consist of an upper and a lower index.
\end{lemma}
\begin{proof}
By our simplifying assumption we may write:
\begin{equation}
\label{powerSeriesForPiS}
{\tilde \Pi}^s(y)^a = y^a + A^a_{\alpha \beta} y^{\alpha} y^\beta + B^a_{\alpha \beta \gamma} y^\alpha y^\beta y^\gamma + O(|y|^4)
\end{equation}
where $A^a_{\alpha \beta}$ is symmetric in $\alpha$ and $\beta$ and $B$ is symmetric in $\alpha$, $\beta$ and $\gamma$. The Taylor series expansion
for $\phi$ now allows us to compute the components of $(y-\phi(\tilde{\Pi}^s(y)))$.
\begin{equation*}
\begin{split}
(y - \phi({\tilde \Pi}^s(y)))^a &= y^a 
-D^a_\alpha
( y^\alpha + A^\alpha_{\beta \gamma} y^\beta y^\gamma + B^\alpha_{\beta \gamma \delta} y^\beta y^\gamma y^\delta) \\
&\quad
-\frac{1}{2}\frac{\partial^2 \phi^a}{\partial x^\alpha \partial x^\beta} 
  ( y^\alpha + A^\alpha_{\delta \epsilon} y^\delta y^\epsilon )( y^\beta + A^\beta_{\zeta \eta} y^\zeta y^\eta ) + O(|y|^4) \\
&= y^a - D^a_\alpha y^\alpha \\
&\quad - D^a_\alpha A^\alpha_{\beta \gamma} y^\beta y^\gamma
-\frac{1}{2} \frac{\partial^2 \phi^a}{\partial x^\alpha x^\beta} y^\alpha y^\beta  \\
&\quad - \frac{\partial^2 \phi^a}{\partial x^\alpha \partial x^\beta} A^\alpha_{\delta \epsilon} y^\delta y^\epsilon y^\beta - D^a_\alpha B^\alpha_{\beta \gamma \delta} y^\beta y^\gamma y^\delta + O(|y|^4)
\end{split}
\end{equation*}
We take the partial derivative of this with respect to $A^p_{qr}$ to get:
\[
\frac{\partial}{\partial A^p_{qr}} (y - \phi({\tilde \Pi}^s(y)))^a
= - D^a_p y^q y^r - \frac{\partial^2 \phi^a}{\partial x^p \partial x^\beta} y^q y^r y^\beta + O(|y|^4).
\]

Because of the distance minimizing property of $\Pi^s$ we know that for all $p$, $q$, $r$ 
and sufficiently small $y$ we have:
\[ \frac{\partial}{\partial A^p_{qr}}|(y - \phi({\tilde \Pi}^s(y))|^2 = 0 \]
The left hand side of this expression is equal to:
\[ 2 \left(\frac{\partial}{\partial A^p_{qr}} (y - \phi({\tilde \Pi}^s(y)))^a \right) (y - \phi({\tilde \Pi}^s(y)))^a. \]
We have written down explicit expressions for each term in this product. This enables us to write
down the fourth order terms of $\frac{\partial}{\partial A^p_{qr}}|(y - \phi({\tilde \Pi}^s(y))|^2$. They
are given by:
\begin{multline*}
2 D^a_p A^a_{\alpha \beta} y^\alpha y^\beta y^q y^r
+ D^a_p \frac{\partial^2 \phi^a}{\partial x^\alpha x^\beta} y^\alpha y^\beta y^q y^r
- 2 \frac{\partial^2 \phi^a}{\partial x^p \partial x^\beta}(y^a - D^a_\gamma y^\gamma) y^\beta y^q y^r \\
= \left( 2 A^p_{\alpha \beta} y^\alpha y^\beta
+ \frac{\partial^2 \phi^p}{\partial x^\alpha \partial x^\beta} y^\alpha y^\beta
- 2 \frac{\partial^2 \phi^a}{\partial x^p \partial x^\beta}(y^a - D^a_\delta y^\delta) y^\beta \right) y^q y ^r \\
= \left( 2 A^p_{\alpha \beta} y^\alpha y^\beta
+ \frac{\partial^2 \phi^p}{\partial x^\alpha \partial x^\beta} y^\alpha y^\beta
- 2 \frac{\partial^2 \phi^a}{\partial x^p \partial x^\beta}(\perp)_{\alpha a} y^\alpha y^\beta \right) y^q y ^r
\end{multline*}
We know that this must vanish for all sufficiently small $y$. We deduce that
\[
2 A^p_{\alpha \beta} y^\alpha y^\beta
+ \frac{\partial^2 \phi^p}{\partial x^\alpha \partial x^\beta} y^\alpha y^\beta
- 2 \frac{\partial^2 \phi^a}{\partial x^p \partial x^\beta}(\perp)_{\alpha a} y^\alpha y^\beta = 0.
\]
for all sufficiently small $y$. This gives us an expression for $A^p_{\alpha \beta} y^\alpha y^\beta$ which combines with equation \eqref{powerSeriesForPiS} to prove the result.
\end{proof}

We now use the lemma coupled with some coordinate transformations to compute a second order expression
for the metric projection  in the general case.
\begin{proposition}
Let $g^\phi_{\perp}$ be the symmetric two form on $\R^r$ defined by:
\[ g^\phi_\perp(X+X^\perp, Y+Y^\perp) = g(X^\perp, Y^\perp) \qquad X,Y \in \Image \phi \hbox{ and } X^\perp,Y^\perp \in (\Image \phi)^\perp
\]
where $g$ is the Euclidean metric on $\R^r$. Define coordinates $\tilde y$ centered
on $0 \in R^r$ by
$\tilde{y}^a=y^a-y_0^a$.
Then to second order the metric projection  is given by
\[
\begin{split}
{\tilde \Pi}^s(y)^a &= x_0 + 
\Pi^a_\alpha  \tilde{y}^\alpha
- \frac{1}{2} \frac{ \partial^2 \phi^\gamma}{\partial x^\alpha \partial  x^\beta}
\Pi^a_\gamma
\Pi^\alpha_\delta  
\Pi^\beta_\epsilon
\tilde{y}^\delta  \tilde{y}^\epsilon \\
&\quad
+ \frac{ \partial^2 \phi^\gamma}{ \partial x^\alpha \partial x^\beta}
\Pi^\beta_\delta
h^{a \alpha} \tilde{y}^\gamma \tilde{y}^\delta
-  \frac{ \partial^2 \phi^\gamma}{ \partial x^\alpha \partial x^\beta}
\Pi^\beta_\epsilon
\Pi^\eta_\gamma
\Pi^\zeta_\delta
h_{\eta \zeta} h^{a \alpha}
\tilde{y}^\delta \tilde{y}^\epsilon \\
&\quad
+ O(|y|^3)
\end{split}
\]
and where $\Pi$ is given by equation \ref{piFirstOrder}.
\end{proposition}
\begin{proof}
We assume without loss of generality that $x_0=0$ and $y_0=\phi(x_0)=0$.

We can find a coordinate transformation $J$ of $\R^n$ which maps an orthonormal basis
of $\R^n$ to the standard basis vectors. We take $x$ to be our original coordinates and
$X$ to be the coordinates obtained by applying $J^{-1}$. So we have:
\[ x^a = J^a_b X^b \]
To satisfy our requirements $J$ must satisfy:
\[ h_{\alpha \beta} J^\alpha_a J^\beta_b = \delta_{a b} \]
Equivalently:
\[ h_{a b} = (J^{-1})^\alpha_a (J^{-1})^\alpha_b. \]
So any pseudo square root of $h_{a b}$ will give an appropriate choice for $J^{-1}$. Taking
the matrix inverse of the above expression we have:
\begin{equation}
h^{a b} = J^a_\alpha J^b_\alpha
\label{pseudoSquareRoot}
\end{equation}
We can now find an orthogonal transformation $T$ of $\R^r$ mapping $\Image \phi*$
to $\R^n \subseteq \R^r$. Hence
$\Phi=T\circ\phi\circ J$
satisfies $(\Phi_*)^a_b = D^a_b$.
We will write $x$ for the original coordinates on $\R^r$ and
define transformed coordinates $X$ by:
\[
X^a = T^a_\alpha x^\alpha.
\]
Let us write $\Pi^\prime$ for the metric projection  associated with the map $\Phi$. The various maps we have just defined are summarized in the commutative diagram below:

\begin{center}
\includegraphics{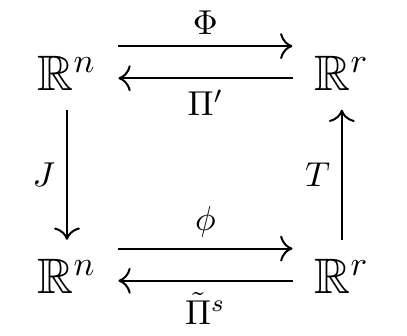}
\end{center}

From Lemma
\ref{secondOrderProjectionLemma} we have:
\begin{align*}
\Pi^\prime(Y)^a &= \frac{\partial \Phi^a}{\partial X^\alpha} Y^\alpha
- \frac{1}{2} \frac{ \partial^2 \Phi^a}{ \partial X^\alpha \partial  X^\beta} Y^\alpha Y^\beta \\
&\quad + \frac{ \partial^2 \Phi^\gamma}{ \partial X^a \partial X^\beta} (\perp)_{\gamma \alpha} Y^\alpha Y^\beta
+ O(|Y|^3) \\
&= \frac{\partial (T^{-1}\circ\Phi)^\beta}{\partial X^\alpha} T^a_\beta Y^\alpha
- \frac{1}{2} \frac{ \partial^2 (T^{-1}\circ\Phi)^\gamma}{\partial X^\alpha \partial X^\beta}T^a_\gamma  Y^\alpha Y^\beta \\
&\quad + \frac{ \partial^2 (T^{-1}\circ\Phi)^\delta}{ \partial X^a \partial  X^\beta} T^\gamma_\delta (\perp)_{\gamma \alpha} Y^\alpha Y^\beta
+ O(|Y|^3). \\
\end{align*}
Hence
\begin{align*}
\Pi^\prime(T y)^a &= \frac{\partial (T^{-1}\circ\Phi)^\beta}{\partial X^\alpha} T^a_\beta T^\alpha_\gamma y^\gamma
- \frac{1}{2} \frac{ \partial^2 (T^{-1}\circ\Phi)^\gamma}{\partial X^\alpha \partial  X^\beta}T^a_\gamma T^\alpha_\zeta T^\beta_\eta  y^\zeta y^\eta \\
&\quad + \frac{ \partial^2 (T^{-1}\circ\Phi)^\delta}{ \partial X^a \partial X^\beta} T^\gamma_\delta (\perp)_{\gamma \alpha} T^\alpha_\zeta T^\eta_\epsilon y^\zeta y^\eta
+ O(|y|^3) \\
&= \frac{\partial (T^{-1}\circ\Phi)^\beta}{\partial x^\gamma} J^\gamma_\alpha T^a_\beta T^\alpha_\gamma y^\gamma
- \frac{1}{2} \frac{ \partial^2 (T^{-1}\circ\Phi)^\gamma}{\partial x^\iota \partial  x^\kappa}
J^\iota_\alpha J^\kappa_\beta
T^a_\gamma T^\alpha_\zeta T^\beta_\eta  y^\zeta y^\eta \\
&\quad + \frac{ \partial^2 (T^{-1}\circ\Phi)^\delta}{ \partial x^\iota \partial x^\kappa}
J^\iota_a J^\kappa_\beta
T^\gamma_\delta (\perp)_{\gamma \alpha} T^\alpha_\zeta T^\beta_\eta y^\zeta y^\eta
+ O(|y|^3). \\
\end{align*}
We deduce that:
\begin{align*}
{\tilde \Pi}^s(y)^a &= ((J \circ \pi^\prime \circ T)(y))^a \\
&= \frac{\partial \phi^\beta}{\partial x^\delta}
J^a_\lambda
J^\delta_\alpha T^\lambda_\beta T^\alpha_\gamma y^\gamma
- \frac{1}{2} \frac{ \partial^2 \phi^\gamma}{\partial x^\iota \partial  x^\kappa}
J^a_\lambda
J^\iota_\alpha J^\kappa_\beta
T^\lambda_\gamma T^\alpha_\zeta T^\beta_\eta  y^\zeta y^\eta \\
&\quad + \frac{ \partial^2 \phi^\delta}{ \partial x^\iota \partial  x^\kappa}
J^a_\lambda
J^\iota_\lambda J^\kappa_\beta
T^\gamma_\delta (\perp)_{\gamma \alpha} T^\alpha_\zeta T^\beta_\eta y^\zeta y^\eta
+ O(|y|^3).
\end{align*}
We now note that:
\[ \Pi^a_b = J^a_\alpha D^\alpha_\beta T^\beta_b =  J^a_\alpha T^\alpha_b. \]
This allows us to simplify our expression for $\Pi^s(y)$ to:
\begin{align*}
{\tilde \Pi}^s(y)^a &= 
\frac{\partial \phi^\beta}{\partial x^\delta}
\Pi^a_\beta
\Pi^\delta_\gamma y^\gamma
- \frac{1}{2} \frac{ \partial^2 \phi^\gamma}{\partial x^\iota \partial  x^\kappa}
\Pi^a_\gamma
\Pi^\iota_\zeta  
\Pi^\kappa_\eta 
y^\zeta y^\eta \\
&\quad + \frac{ \partial^2 \phi^\delta}{ \partial x^\iota \partial x^\kappa}
\Pi^\kappa_\eta
J^a_\lambda
J^\iota_\lambda
T^\gamma_\delta (\perp)_{\gamma \alpha} T^\alpha_\zeta y^\zeta y^\eta
+ O(|y|^3).
\end{align*}
The tensor $\perp_{a b}$ is equal to $(g^\Phi_{\perp})_{ab}$. 
So since $T$ is an isometry, we may write
\[
(\perp)_{\alpha \beta} T^\alpha_a T^\beta_b = (g^\phi_{\perp})_{ab}.
\]
Using this together with equation \eqref{pseudoSquareRoot} we may write:
\begin{align*}
{\tilde \Pi}^s(y)^a &= 
\frac{\partial \phi^\beta}{\partial x^\delta}
\Pi^a_\beta
\Pi^\delta_\gamma y^\gamma
- \frac{1}{2} \frac{ \partial^2 \phi^\gamma}{\partial x^\iota \partial  x^\kappa}
\Pi^a_\gamma
\Pi^\iota_\zeta  
\Pi^\kappa_\eta 
y^\zeta y^\eta \\
&+ \frac{ \partial^2 \phi^\delta}{ \partial x^\iota \partial x^\kappa}
\Pi^\kappa_\eta
h^{a \iota}
(g^\phi_{\perp})_{\delta \zeta} y^\zeta y^\eta
+ O(|y|^3)
\end{align*}
The first term can be simplified by repeated applications of equations \eqref{metricTensor}
and \eqref{piFirstOrder}:
\begin{align*}
\frac{\partial \phi^\beta}{\partial x^\delta}
\Pi^a_\beta
\Pi^\delta_\gamma y^\gamma &= 
\frac{\partial \phi^\beta}{\partial x^\delta}
\frac{\partial \phi^\beta}{\partial x^\gamma} h^{a \gamma}
\frac{\partial \phi^\gamma}{\partial x^\epsilon} h^{\epsilon \delta} y^\gamma \\
& = h_{\delta \gamma} h^{\alpha \gamma} \frac{\partial \phi^\gamma}{\partial x^\epsilon} h^{\epsilon \delta} y^\gamma
=  \frac{\partial \phi^\gamma}{\partial x^\epsilon} h^{\epsilon a} y^\gamma
= \Pi^a_\beta y^\beta.
\end{align*}
It is a tautology that the first order term is given by $\Pi$, nevertheless this calculation is
a reassuring check on our working.
Renaming the dummy variables we now have that:
\[
\begin{split}
{\tilde \Pi}^s(y)^a &= x_0 + 
\Pi^a_\alpha  y^\alpha
- \frac{1}{2} \frac{ \partial^2 \phi^\gamma}{\partial x^\alpha \partial  x^\beta}
\Pi^a_\gamma
\Pi^\alpha_\delta  
\Pi^\beta_\epsilon
{ y}^\delta  y^\epsilon \\
&\quad
+ \frac{ \partial^2 \phi^\gamma}{ \partial x^\alpha \partial x^\beta}
\Pi^\beta_\epsilon
h^{a \alpha}
(g^\phi_{\perp})_{\gamma \delta} y^\delta y^\epsilon
+ O(|y|^3)
\end{split}
\]
We would like a formula that can be computed efficiently when $n \ll r$, so we wish to eliminate the term $g^\phi_{\perp}$. By splitting vectors $V$ and $W$ in $\R^r$ into components in $\Image \phi_*$ and its orthogonal complement, we see that the Euclidean metric on $\R^r$ satisfies the decomposition:
\[ g_{ab} V^a W^b = (g^\phi_{\perp})_{ab} V^a W^b + h_{ab} \Pi^a_\alpha V^\alpha \Pi^b_\alpha  W^\beta.  \]
Using this formula we obtain:
\[
\begin{split}
{\tilde \Pi}^s(y)^a &= x_0 + 
\Pi^a_\alpha  y^\alpha
- \frac{1}{2} \frac{ \partial^2 \phi^\gamma}{\partial x^\alpha \partial  x^\beta}
\Pi^a_\gamma
\Pi^\alpha_\delta  
\Pi^\beta_\epsilon
{ y}^\delta  y^\epsilon \\
&\quad
+ \frac{ \partial^2 \phi^\gamma}{ \partial x^\alpha \partial x^\beta}
\Pi^\beta_\delta
h^{a \alpha} y^\gamma y^\delta
-  \frac{ \partial^2 \phi^\gamma}{ \partial x^\alpha \partial x^\beta}
\Pi^\beta_\epsilon
\Pi^\eta_\gamma
\Pi^\zeta_\delta
h_{\eta \zeta} h^{a \alpha}
y^\delta y^\epsilon \\
&\quad
+ O(|y|^3)
\end{split}
\]
\end{proof}

We can immediately conclude:
\begin{theorem}[It\^o-jet projection in coordinates]\label{explicitVector}
Let $\phi:\R^n \to \R^r$ be an embedding with $\phi(x_0)=y_0$
then the It\^o-jet projection of the SDE:
\[
\ed y = a \, \ed t + b_\alpha \, \ed W^\alpha_t,     \qquad y_0
\]
is
\[
\ed x = A \, \ed t + B_\alpha \, \ed W^\alpha_t,     \qquad x_0
\]
where:
\[
B^i_\alpha = \Pi^i_\beta b^\beta_\alpha
\]
and:
\[
\begin{split}
A^i &= \Pi^i_\alpha a^\alpha + \\
&\quad
\left( -\frac{1}{2} \frac{ \partial^2 \phi^\gamma}{\partial x^\alpha \partial  x^\beta}
\Pi^i_\gamma
\Pi^\alpha_\delta  
\Pi^\beta_\epsilon
\right. \\
&\quad
\left.
+ \frac{ \partial^2 \phi^\epsilon}{ \partial x^\alpha \partial x^\beta}
\Pi^\beta_\delta
h^{i \alpha} 
-  \frac{ \partial^2 \phi^\gamma}{ \partial x^\alpha \partial x^\beta}
\Pi^\beta_\epsilon
\Pi^\eta_\gamma
\Pi^\zeta_\delta
h_{\eta \zeta} h^{i \alpha}
\right) \\
&\quad \times b^\delta_\kappa b^\epsilon_\iota [W^\kappa, W^\iota]_t.
\end{split}
\]
$\Pi$ is given by \eqref{piFirstOrder}. $h_{ab}$ is given by \eqref{metricTensor}. $h^{ab}$ is the inverse of $h_{ab}$.
\end{theorem}

It is reassuring to check that this formula gives the
same result as we found in \Cref{section:lowDimensionalExample} for
projection of a particular SDE onto a circle where the projection map
was known exactly. In fact, we can find an explicit expression for the It\^o-jet projection of any bivariate SDE driven by a single Brownian motion on the plane on the unit circle. 

\begin{example}[It\^o-jet projection of a bivariate SDE on the unit circle]
Suppose that our diffusion process in $\R^r=\R^2$, driven by a one-dimensional Brownian motion $W=W^1$, is 
\[ dX = a_1(X,Y) dt + b^1_1(X,Y) dW^1, \ X_0 \]
\[ dY = a_2(X,Y) dt + b^2_1(X,Y) dW^1, \ Y_0\]
and suppose we wish to approximate this process in the unit circle. 
If we define $\theta = \arctan(Y_t/X_t)$, and compute $d \theta_t$ via It\^o's formula, this won't be in general a closed SDE for $\theta$, contrary to the special example of the cross diffusion above. To obtain a closed SDE in $\theta$ we have to project. One can check that for the one-dimensional manifold given by the unit circle, expressed as
\[ M= \{ (\cos(\theta), \sin(\theta)), \theta \in [0, 2\pi) \} \]
with coordinates $Y=\theta$ in $\R^n=\R^1$, one has
\[ h=1, \ h^{-1}=1, \ \Pi_1 = - \sin(\theta), \ \Pi_2 = \cos(\theta), \ \ \partial^2_{\theta} \phi^1 = - \cos(\theta), \ \ 
\partial^2_{\theta} \phi^2 = - \sin(\theta), \]
which allows us to apply \Cref{explicitVector} to this system. We obtain (coefficients $a$ and $b$ are computed in $X=\cos(\theta), Y=\sin(\theta)$)
\[ A(\theta) = - a_1 \sin(\theta) + a_2 \cos(\theta) + \frac{1}{2} \sin(2 \theta) ((b^1_1)^2-(b^2_1)^2) - \cos(2 \theta) b^1_1 b^2_1, \]
\[ B(\theta) = - \sin(\theta) b^1_1 + \cos(\theta) b^2_1.\]
In the special case of $a_1=a_2=0$ and $b^1_1 = \sigma \sin(\theta), b^2_1 = \sigma \cos(\theta)$ this confirms our previous calculations for the cross-diffusion example.
\end{example}

\section{Application of the Projection to Non-linear Filtering}
\label{applicationSection}
As a fundamental application of our new projection methods we consider an area from signal processing, stochastic filtering. This extends our previous work in \cite{armstrongbrigomcss}.

In stochastic filtering one has a signal $X$ that evolves according to a SDE, and observes a process $Y$ which is a function of this signal plus noise. This is standard notation, but these $X$ and $Y$ are not to be confused with the processes we used earlier in the paper, in that they are not the $\R^r$ process to be approximated and its $\R^n$ approximation.

The filtering problem consists in estimating the signal $X$ given the present and past observations $Y$. If $t$ is the current time, the solution of the filtering problem is the probability density of the state $X_t$ conditional on the observations from time 0 to time $t$, call it $p_t$. The density $p_t$ follows the Kushner-Stratonovich (or Zakai) stochastic partial differential equation (SPDE) that, under some technical assumptions, can be seen as a stochastic differential equation in the infinite dimensional $L^2$ space of square roots of densities (Hellinger metric) or of densities themselves (direct $L^2$ metric).

The process we wish to approximate on a low dimensional manifold is $p_t$, which represents the $X_t$ of our earlier sections. The $\R^r$ space of our earlier sections is the $L^2$ infinite dimensional space, while the submanifold $M$ is a finite dimensional family of probability densities parametrized by $\theta$, acting as coordinates: $\{p(\cdot,\theta), \ \theta \in \Theta \subset \R^n\}$. $\theta_t$ plays the role of what we were calling $Y_t$ earlier in the paper. We aim at finding a SDE for $\theta$ such that $p(\cdot,\theta_t)$ approximates $p_t(\cdot)$ in an optimal way. Note that in the previous part of the paper we had a dimensionality reduction from $r$ to $n$, whereas now we go from infinite dimensional $p_t$ to $n$-dimensional $\theta_t$.

One may be concerned about taking our finite dimensional results and
applying them in an infinite dimensional setting. However, we have stated our
results in terms of approximating one Ito--Taylor series of a given order with
another Ito--Taylor series. This allows us to avoid the analytical issues that
might conceivably arise in considering the convergence of these series. Therefore our results generalize straightforwardly to the Hilbert space setting.
As an example, the minimization argument used to prove \Cref{theoremStrongVector} relies only on properties of the linear projection operator that remain true in a Hilbert space setting.

In addition the explicit calculation of \Cref{section:secondOrderProjection} can be generalized unproblematically to the case of a finite dimensional manifold embedded in a Hilbert space. To see this simply note that the vector space spanned by the first two derivatives of the map $\phi$ at $p$ gives a finite dimensional space $V$ and so one can simply apply the result for embedding into the space $V$.

The point where complexities might conceivably arise in the infinite dimensional
setting is in the generalizations of \Cref{proposition:strongConvergence}
and \Cref{proposition:weakConvergence}. Folk wisdom suggests that such results
can be generalized to Hilbert spaces without difficulty, so we will not attempt to prove that here.

\subsection{The Kushner Stratonovich equation}
\label{filteringSubsection}

We suppose that the state $X_t \in \R^n$ of a system
evolves according to the equation:
\[ \ed X_t = f(X_t,t) \, \ed t + \sigma(X_t,t) \, \ed W_t \]
where $f$ and $\sigma$ are smooth $\R^n$ valued functions
and $W_t$ is a Brownian motion. One typically adds growth conditions to ensure a global existence and uniqueness result for the signal equation, see for example \cite{armstrongbrigomcss} and references therein for the details.

We suppose that an associated process, the observation process, $Y_t \in \R^d$
evolves according to the equation:
\[ \ed Y_t = b(X_t,t) \, \ed t + \ed V_t \]
where $b$ is a smooth $\R^d$ valued function and $V_t$ is a Brownian motion independent of $W_t$.
Note that the filtering problem is often formulated with an additional constant in terms
of the observation noise. For simplicity we have assumed that the system is scaled so that this can be omitted.

The filtering problem is to compute the conditional
distribution of $X_t$ given a prior distribution for $X_0$
and the values of $Y$ for all times up to and including~$t$.

Subject to various bounds on the growth of the
coefficients of this equation, the assumption that
the distribution has a density $p_t$ and suitable
bounds on the growth of $p_t$ one can show that
$p_t$ satisfies the Kushner--Stratonovich SPDE:
\begin{equation}
\ed p = {\cal L^*}  p \  \ed t 
+ p[b - E_p(b)]^T [ \ed Y - E_p(b) \ed t]
\label{KSIto}
\end{equation}
where $E_p$ denotes the expectation with respect to
the density $p$, 

$E_p[f] = \int f(x) p(x) dx$, and the forward diffusion operator ${\cal L}^*_t$ is defined by:  
\begin{equation}
{\cal L}_t^* \phi = - \frac{\partial}{\partial x^i} [ f_i(x,t) \phi ] + \frac{1}{2} \frac{\partial^2}{\partial x^i \partial x^j} [a_{ij}(x,t)\phi]
\end{equation}
where $a=\sigma \sigma^T$.
Note that we are using the Einstein summation convention in this expression. 

In the event that the coefficient functions $f$
and $b$ are all linear and $\sigma$ is a deterministic function of time 
one can show that
so long as the prior distribution for $X$ is Gaussian, or deterministic, 
the density $p$ will be Gaussian at all subsequent times.
This allows one to reduce the infinite dimensional 
equation \eqref{KSIto} to a finite dimensional stochastic differential equation for the mean and covariance matrix of this normal distribution. This finite dimensional
problem solution is known as the Kalman filter.

For more general coefficient functions, however, equation
\eqref{KSIto} cannot be reduced to a finite dimensional
problem \cite{hazewinkel}. Instead one might seek approximate solutions of \eqref{KSIto} that belong to some given statistical family of densities. This is
a very general setup and includes, for example,
approximating the density using piecewise linear
functions to derive a finite difference approximation
or approximating the density with Hermite polynomials
to derive a spectral method. Other examples include exponential families (considered in \cite{brigo2,brigo1}) and mixture families (considered in \cite{armstrongBrigo,armstrongbrigomcss}).

Our projection theory tells us how one can find good approximations 
on a given statistical family with respect to a given
metric on the space of distributions. We illustrate
this by writing down the It\^o-vector and It\^o-jet projection of
\eqref{KSIto} for the $L^2$ and Hellinger metrics
onto a general manifold\footnote{Note that it is also possible to consider projecting the Zakai equation. However, as explained in 
\cite{armstrongbrigomcss}, one expects that projecting the Kushner--Stratonovich will lead to
smaller error terms.}. 

We will then examine some numerical results regarding 
the very specific case
of seeking approximate solutions using Gaussian distributions. The idea of approximating the solution
to the filtering problem using a Gaussian
distribution has been considered by numerous authors
who have derived variously, the extended Kalman filter
\cite{pardoux}, assumed density filters \cite{kushner} and Stratonovich projection filters \cite{brigo1}. Some of these are related, for example the assumed density filters and Stratonovich projection filters in Hellinger metrics for Gaussian (and more generally exponential) families coincide \cite{brigo2}.  Using our new projection methods,  we will be able to derive projection filters which
 outperform all these other filters
(assuming performance is measured over small time intervals using the appropriate Hilbert space metric).

We note that \eqref{KSIto} is an infinite dimensional SDE
driven by a continuous semi-martingale. The definitions and results
given in \Cref{section:taylorSeries} were only stated in the
finite dimensional case for SDEs driven by Brownian motion.
The definition of It\^o--Taylor series can be generalized straightforwardly
to this situation and hence the definition of the It\^o projections can
be applied in this context also.

More generally, for the the geometry of infinite dimensional filtering problems based on $L^2$ or Orlicz charts and for the related differential geometric approach to statistics with recent advances we refer for example to \cite{pistoneannals,newton1,newton2,girolami,brigo2,armstrongbrigomcss,brigopistone,brigopistone2}.

\subsection{It\^o-vector projections}
\label{vectorCalculations}
\subsubsection{The It\^o-vector projection filter in the $L^2$ direct metric}

Let us suppose that the density $p$ lies in $L^2$ and so
we can use the $L^2$ norm to measure the accuracy of
an approximate solution to equation \eqref{KSIto}. For a discussion on conditions under which a unnormalized version of $p$ is in $L^2$ (Zakai Equation) see for example \cite{ahmedbook}. 

We wish to consider an $m$-dimensional family
of distributions $p$ parameterized by $m$ real valued parameters $\theta^1$, $\theta^2$, $\ldots$, $\theta^m$. For example we will consider the 2 dimensional Gaussian family:
\begin{equation}
 p(x) = \frac{1}{(\theta^2) \sqrt{2 \pi}}\exp\left(-\frac{(x-(\theta^1))^2}{2 (\theta^2)^2}\right). 
 \label{gaussianFamily}
\end{equation}
Note that we have chosen to follow differential
geometry convention and use upper indices for
the coordinate functions $\theta^i$ so we have 
been careful to distinguish powers from indices using
brackets.

More formally, an $m$-dimensional family is given
by a smooth embedding $\phi:\R^m \rightarrow L^2(\R^n)$.
The tangent vectors $\phi_*\frac{\partial}{\partial \theta^i} \in L^2(\R^n)$ are simply the partial derivatives
\[
\frac{\partial p}{\partial \theta^i}.
\]

Let us write:
\[ g_{ij}=\int_{\R} \frac{\partial p}{\partial \theta^i}
\frac{\partial p}{\partial \theta^j} \, \ed x.\]
This defines the induced metric tensor on the manifold
$\phi(\R^m)$. We will write $g^{ij}$ for the inverse of the matrix $g_{ij}$. The projection operator $\Pi_{\phi(\theta)}$
is then given by
\begin{equation*}
\begin{split}
\Pi_{\phi(\theta)} (v) &= 
\sum_{i,j=1}^m g^{ij} \left\langle v, \phi_* \frac{\partial}{\partial \theta^i} \right\rangle_{L^2}
\phi_*\frac{\partial}{\partial \theta^j} \\
&= \sum_{i,j=1}^m g^{ij} \left( \int_{\R^n} v(x) \frac{\partial p}{\partial \theta^i} \, \ed x \right)
\phi_*\frac{\partial}{\partial \theta^j} \, .
\end{split}
\end{equation*}
Thus
\[
\phi_*^{-1} \Pi_{\phi(\theta)} (v)
= \sum_{i,j=1}^m g^{ij} \left( \int_{\R^n} v(x) \frac{\partial p}{\partial \theta^i} \, \ed x  \right)
\frac{\partial}{\partial \theta^j} \, .
\]
We can now write down the It\^o-vector projection of \eqref{KSIto} with respect to the $L^2$ metric. It is:
\[ \ed \theta^i = A^i \, \ed t + B^i \, \ed Y_t \]
where:
\[ B^i =
\sum_{j=1}^m g^{ij} \left( \int_{\R}
(p(b-E_{p(\theta)}(b)))^T
\frac{\partial p}{\partial \theta^j} \, \ed x \right)
\]
and
\begin{small}
\[ A^i =
\sum_{j=1}^m g^{ij} \left( \int_{\R^n}
\left(
{\cal L}^*p - p(b-E_{p(\theta)}(b))^T E_{p(\theta)}(b)
- \frac{1}{2} \sum_{k=1}^m
\frac{\partial^2 p}{\partial \theta^j \partial \theta^k} B^k
\right)
\frac{\partial p}{\partial \theta^j} \, \ed x \right).
\]
\end{small}

\begin{example}[It\^o-vector projection filter for cubic sensor in direct metric]
\label{cubicExample}
Consider as a test case the
$1$-dimensional problem with $f(x,t)=0$, $\sigma(x,t)=1$
and $b(x,t)=x + \epsilon x^3$ for some small constant $\epsilon$. This problem is a perturbation of a linear filter so one might expect that a Gaussian approximation will perform reasonably well at least for small times. Thus
we will use the 2 dimensional manifold of Gaussian distributions given in equation \eqref{gaussianFamily}.

We first calculate the metric tensor $g_{ij}$ which
is diagonal in this case:
\begin{equation*}
g_{ij} = \frac{1}{4 \sqrt{\pi } (\theta^2)^3} \left(
\begin{array}{cc}
 1 & 0 \\
 0 & \frac{3}{2} \\
\end{array}
\right) \, .
\end{equation*}
This is easily inverted to compute $g^{ij}$. We
compute the expectation $E_p(b)$:
\begin{equation*}
E_p(b) = 
\frac{\epsilon \left(\sqrt{2 \pi } (\theta^1)^3 (\theta^2)+3 \sqrt{2 \pi } (\theta^1) (\theta^2)^3\right)}{\sqrt{2 \pi } (\theta^2)}+(\theta^1).
\end{equation*}
One can now see that computing the projection equation
will simply involve integrating a number of terms
of the form a polynomial in $x$ times a Gaussian.
The end result is:
\begin{scriptsize}
\itoLTwoProjectedEquation
\end{scriptsize}
\end{example}

\subsubsection{The It\^o-vector projection filter in the Hellinger metric}

The Hellinger metric is a metric on probability
measures. In the case of two probability
density functions $p(x)$ and $q(x)$ on $\R^n$, that now need only be in $L^1$, the Hellinger distance is
given by the square root of:
\[ 
\frac{1}{2} \int (\sqrt{p(x)} - \sqrt{q(x)})^2 \, \ed x .
\]
In other words, up to the constant factor of $\frac{1}{2}$
the Hellinger metric corresponds to the $L^2$ norm on
the square root of the density function rather than on the density itself (as in the previous subsection). The Hellinger metric has the important advantage of making the metric independent of the particular background density that is used to express measures as densities. The $L^2$ direct distance introduced earlier does not satisfy this background independence.

Now, to compute the It\^o-vector projection with respect to
the Hellinger metric we first want to write down
an It\^o equation for the evolution on $\sqrt{p}$.

Applying It\^o's lemma to equation \eqref{KSIto} we
formally obtain:
\begin{equation*}
\begin{split}
\ed \sqrt{p}
&= \left( 
\frac{ {\cal L}^* p - p(b - E_p(b))^T E_p(b)}{ 2 \sqrt{p}}
- \frac{p^2 (b - E_p(b))^T(b - E_p(b))}{ 8 p \sqrt{p}}
\right) \, \ed t \\
&{} + \left(
\frac{p(b - E_p(b))^T }{ 2 \sqrt{p}}
\right) \ed Y_t. \\
&= \left( 
\frac{ {\cal L}^* p}{ 2 \sqrt{p}}
- \frac{1}{8}\sqrt{p} (b - E_p(b))^T( b + 3 E_p(b))
\right) \, \ed t \\
&{} + \left(
\frac{1}{2} \sqrt{p}(b - E_p(b))^T
\right) \ed Y_t.
\end{split}
\end{equation*}

A family of distributions now corresponds to an embedding
$\phi$ from $\R^m$ to $L^2(\R^n)$ but now $p=\phi(\theta)^2$. The tangent space is spanned by
the vectors:
\[ \phi_* \frac{\partial}{\partial \theta^i}
 = \frac{ \partial \sqrt{p} }{\partial \theta^i}. \]
We define a metric on the tangent space by:

\[
h_{ij} = \int_{\R^n} \frac{ \partial \sqrt{p} }{\partial \theta^i} \frac{ \partial \sqrt{p} }{\partial \theta^j} \, \ed x.
\]
We write $h^{ij}$ for the inverse matrix of $h_{ij}$. The projection operator with respect to the Hellinger metric is:
\begin{equation*}
\begin{split}
\Pi_{\phi(\theta)} (v) 
&= \sum_{i,j=1}^m h^{ij} \left( \int_{\R^n} v(x) \frac{\partial \sqrt{p}}{\partial \theta^i} \, \ed x \right)
\phi_*\frac{\partial}{\partial \theta^j} .
\end{split}
\end{equation*}
We can now write down the It\^o-vector projection of \eqref{KSIto} with respect to the Hellinger metric. It is: 
\[ \ed \theta^i = A^i \, \ed t + B^i \, \ed Y_t \]
where:
\[ B^i =
\sum_{j=1}^m h^{ij} \left( \int_{\R}
\frac{1}{2} \sqrt{p}(b - E_{p(\theta)}(b))^T
\frac{\partial \sqrt{p}}{\partial \theta^j} \, \ed x. \right)
\]
and
\[
\begin{split}
A^i &=
\sum_{j=1}^m h^{ij} \left( \int_{\R^n}
\left(
\frac{{\cal L}^* p}{ 2 \sqrt{p}}
- \frac{1}{8}\sqrt{p} (b - E_{p(\theta)}(b))^T( b + 3 E_{p(\theta)}(b)) \right. \right. \\
&\qquad \left. \left.
- \frac{1}{2} \sum_{k=1}^m
\frac{\partial^2 \sqrt{p}}{\partial \theta^j \partial \theta^k} B^k
\right)
\frac{\partial \sqrt{p}}{\partial \theta^j} \, \ed x. \right).
\end{split}
\]

\begin{example}[It\^o-vector projection filter for cubic sensor: Hellinger metric]
We may repeat example \ref{cubicExample} but projecting
using the Hellinger metric.
We first calculate the metric tensor $h_{ij}$ which
is diagonal also in this case:
\begin{equation*}
h_{ij} =  \frac{1}{4 \theta_2^2} \left(
\begin{array}{cc}
 1 & 0 \\
 0 & 2 \\
\end{array}
\right)
\end{equation*}
This is easily inverted to compute $h^{ij}$.
We obtain the following SDEs:
\begin{scriptsize}
\itoHellingerProjectedEquation
\end{scriptsize}
\end{example}

\subsection{It\^o-jet projections}

Using the formulae from \Cref{explicitVector} together with the
formulae and techniques of \Cref{vectorCalculations} we can
explicitly calculate the It\^o-vector projections of the filtering
equation in both the $L^2$ and Hellinger metrics.

To minimize notation, let us concentrate on the $1$-dimensional state space
filtering problem and project using the $L^2$ metric.

We can formally write the filtering equation in the form:
\begin{equation}
\ed p_t = \mu(p_t) \ed t + \Sigma(p_t) \ed W_t
\label{infDimEquation}
\end{equation}
where $p_t$ is an $L^2$ function and
\begin{equation}
\label{filteringCoefficientsInfDim}
\begin{split}
\mu(p)(x)&:= \frac{1}{2} \frac{ \ed^2 (\sigma(x)^2 p(x))}{ \ed x^2}
		   - \frac{\ed (f(x) p(x))}{\ed  x} \\
&{}\quad   -  p(x) \left( b(x) - \int_\R p(t) b(t) \ed t \right) \int_{\R} p(t)b(t) \ed t, \\
\Sigma(p)(x)&:= p(x) \left( b(x) - \int_{\R} p(t) b(t) \ed t \right).
\end{split}
\end{equation}
We now suppose that $p_t$ is parameterized as $p_t(x) = \phi(\theta)(x)$
as in \Cref{vectorCalculations}. Using \Cref{explicitVector}
we can write down the It\^o-jet projection which is an SDE for the components of $\theta$.

To write down the result it will be useful to define functions $\pi^i(\theta)$ by:
\[
\pi^i(\theta) = h^{ij} \frac{\partial \phi}{\partial \theta^j} (\theta).
\]
We will also use 
angle brackets to denote the $L^2$ inner product. With this
understood, the It\^o-jet projection of the filtering
equations in the $L^2$ metric is given by:
\[ 
\ed \theta^i_t = A^i(\theta) \ed t + B^i(\theta) \ed W_t
\]
where we have in turn
\[
B^i(\theta) = \langle \pi^i(\theta), \Sigma(\phi(\theta)) \ed t
\]
and
\begin{equation*}
\begin{split}
A^i(\theta) &= \langle \pi^i(\theta), \mu(\phi(\theta)) \rangle \\
& 
\quad - \frac{1}{2}
\left\langle \frac{\partial^2 \phi}{\partial \theta^\alpha \partial \theta^\beta}(\theta),
	      \pi^{i}(\theta)
\right\rangle
\langle
\Sigma(\theta), \pi^\alpha(\theta)
\rangle
\langle
\Sigma(\theta), \pi^\beta(\theta)
\rangle \\
& 
\quad + 
\left\langle \frac{\partial^2 \phi}{\partial \theta^\alpha \partial \theta^\beta}(\theta), \Sigma(\phi(\theta)) \right\rangle
\langle \pi^\beta(\theta), \Sigma(\phi(\theta)) \rangle
h^{i \alpha}(\theta) \\
& 
\quad -
\left\langle \frac{\partial^2 \phi}{\partial \theta^\alpha \partial \theta^\beta}, \pi^\eta(\theta) \right\rangle
\langle \pi^\beta(\theta), \Sigma(\phi(\theta)) \rangle
\langle \pi^\xi(\theta), \Sigma(\phi(\theta)) \rangle
h_{\eta \xi}(\theta) h^{i \alpha}(\theta)
\end{split}
\end{equation*}

\begin{example}[It\^o-jet projection filter for cubic sensor in direct metric]
For the filtering problem of \ref{cubicExample} the
It\^o-jet projection in the $L^2$ metric is
\begin{scriptsize}
\itoLTwoProjectedEquationB
\end{scriptsize}
\end{example}

The It\^o-jet projection of the filtering equation in the Hellinger metric can be computed in the same way. Indeed we can formally write the filtering equation in the form:
\begin{equation}
\ed q_t = \mu(q_t) \ed t + \Sigma(q_t) \ed W_t
\end{equation}
where $q_t$ is the square root of the density and the coefficients now
satisfy
\begin{equation}
\label{filteringCoefficientsInfDimHellinger}
\begin{split}
\mu(q)(x)&:= \frac{1}{2 q(x)}
\left(\frac{1}{2} \frac{d^2 (\sigma(x)^2 q(x)^2)}{\ed x^2} - \frac{\ed (f(x)q(x)^2)}{\ed x} \right) \\
&\quad
  - 
 \frac{1}{8} q(x) \left( b(x) - \int_\R q(t)^2 b(t) \ed t \right)
\left( b(x) + 3 \int_{\R} q(t)^2 b(t) \ed t \right), \\
\Sigma(q)(x)&:= \frac{1}{2} q(x) \left( b(x)- \int_{\R} q(t)^2 b(t) \ed t \right).
\end{split}
\end{equation}
Thus we can use the same formulae as above to compute the Hellinger projection
except we must use the coefficients from \eqref{filteringCoefficientsInfDimHellinger} rather than those
from \eqref{filteringCoefficientsInfDim}.

\begin{example}[It\^o-jet projection filter for cubic sensor: Hellinger metric]
For the filtering problem of \ref{cubicExample}, the
It\^o-jet projection in the Hellinger metric is
\begin{scriptsize}
\itoHellingerProjectedEquationB
\end{scriptsize}
\end{example}

\subsection{Other Gaussian Approximate Filters}

Many other Gaussian approximate filters have been proposed in the
past. We will briefly review a number of different Gaussian approximate filters that can be found
in the literature and calculate the relevant stochastic differential equations for our
example \ref{cubicExample}. We will then compare the performance of these filters numerically.

\subsubsection{The Stratonovich projection filter} 

Instead of using the It\^o-vector projection, one can use the Stratonovich projection.

\begin{example}[Stratonovich proj.\ filter for cubic sensor: direct metric]
General formulae for performing the Stratonovich $L^2$ projection are given in \cite{armstrongBrigo}.
In the specfic case of example \ref{cubicExample} the resulting It\^o SDEs are:

\begin{scriptsize}
\stratLTwoProjectedEquation
\end{scriptsize}
\end{example}

\begin{example}[Stratonovich proj.\ filter for cubic sensor: Hellinger metric]
General formulae for performing the Stratonovich Hellinger projection are given in \cite{brigo2}. 
In the specfic case of example \ref{cubicExample} the resulting SDEs are:
\begin{scriptsize}
\stratHellingerProjectedEquation
\end{scriptsize}
\end{example}

\subsubsection{The Extended Kalman Filter}

The Extended Kalman Filter (EKF) is a heuristically derived method of
finding approximate solutions to the filtering problem based on the idea
of linearising the problem and then using the solution to the linear problem.
In particular one assumes that the solution can be well approximated by a Gaussian distribution.
For the EKF see \cite{jazwinski,ahmedbook}. A definition and heuristic derivation is given in \cite{crisan} (which is based, in turn, on
the derivation given in \cite{pardoux}).


The EKF can be shown to work well on condition that the initial position of
the signal is approximated well, the non-linearities of $f$ are small, $b$
is injective and the observation noise is small \cite{picard}. Moreover, the EKF is widely
used in practice, see \cite{crisan} for references to applications.

\begin{example}
For the example problem $b(x)=x+\epsilon x^3$ the EKF is:
\begin{footnotesize}
\ekfEquation
\end{footnotesize}
\end{example}

\subsubsection{Assumed density filters}

Assumed density filters (ADFs) provide a finite dimensional method of finding
approximate solutions to the filtering problem. They have been considered in, for example,
\cite{kushner}, \cite{maybeck} and \cite{brigo2}.

The general setup is to consider a statistical family $\pi(\cdot,\eta)$ of 
probability measures parameterized by some coordinates $\eta=(\eta^1, \ldots, \eta^m)$. This parameterization is not arbitrary. It must be chosen in such a way that, for elements of the statistical family, the values of $\eta$ correspond to the expectations of some 
twice differentiable scalar functions $\{c^1, \ldots, c^m\}$ defined on $\R^n$.
\[ \eta^i = E_{\pi(\cdot, \eta)}(c^i) =: E_{\eta_t}(c^i) \]
where for brevity we are using the abbreviation $E_{\eta_t}$ for $E_{\pi(\cdot, \eta)}$.

For example one might take the statistical family of normal distributions parameterized
by its first and second moments $\eta_1$ and $\eta_2$, so $c^1(x)=x$, $c^2(x)=(x)^2$.

Given a statistical family parameterized in this way, we define the It\^o ADF to be:
\[
\ed \eta^i_t = E_{\eta_t}({\cal L}_t c^i) \, \ed t
+ \left( E_{\eta_t}(b_t c^i) - E_{\eta_t}(b_t) \eta^i_t \right)^T
\left( \ed Y_t - E_{\eta_t}( b_t) \, \ed t \right).
\]
This is motivated by the fact that under the conditions used
to derive equation \eqref{KSIto}, we have that the $c_i$-moments of $\pi_t$, the true solution to the filtering problem, satisfy the It\^o equation:
\[
\ed \pi_t(c_i) = \pi_t( {\cal L}_t c^i ) \, \ed t - \frac{1}{2}( \pi_t( b c^i ) - \pi_t(b)\pi_t( c^i ))( \ed Y - \pi_t(b) \, \ed t).
\]
Thus if it were true that the true density was a member of our chosen statistical family then
the It\^o ADF would certainly be satisfied. One just hopes that the It\^o ADF will continue to give a reasonable approximation even though we know that the true density isn't a member of 
the chosen statistical family.

With a similar motivation we define the Stratonovich ADF to be:
\[
\begin{split}
\ed \eta^i_t &= E_{\eta_t}( {\cal L}_t c^i ) \, \ed t -\frac{1}{2}( E_{\eta_t}( |b_t|^2 c^i ) - E_{\eta_t}(|b_t|^2)\eta^i_t) \, \ed t \\
&\quad
+ \left( E_{\eta_t}( b c^i ) - E( b_t ) \eta^i_t \right)^T \, \circ \ed Y_t.
\end{split}
\]
If it were true that the density was a member of our statistical family then the It\^o
ADF and the Stratonovich ADF would be equivalent equations. Since we only expect to be able to approximate the true density with our statistical family, we must expect that the It\^o ADF and Stratonovich ADF are in fact inequivalent equations. Intuitively, we can say that the local moment matching approximation on which the ADF heuristics are based and the It\^o-Stratonovich transformation do not commute. 

The justification just given for ADFs is far from convincing. We are relying on little other than hope that these equations will give good approximations. However, it was shown
in \cite{brigo2} that in fact for exponential families, the Stratonovich projection filter in the Hellinger metric coincides with the Stratonovich ADF, and in \cite{brigosmall} that for the Gaussian case, this filter approaches the optimal filter under small observation noise. 
\cite{brigopistone2} show that the equivalence between projection using the $L^2$ direct structure and the assumed density approximation holds for the prediction step of the filtering algorithm, namely the Kolmogorov equation, when using mixture families.

\begin{example}
If we calculate the It\^o assumed density filter corresponding
to example \ref{cubicExample} and the family of normal distributions, and then
change coordinates to $\theta^1$ and $\theta^2$ as used in the previous examples, 
we obtain the SDEs:
\begin{scriptsize}
\itoADFEquation
\end{scriptsize}
\end{example}
\begin{example}
The family of normal distributions is an exponential family, therefore the
Stratonovich assumed density filter is equivalent to the Stratono\-vich projection
filter in the Hellinger metric.
\end{example}

\subsection{Results}

Our explicit calculations show that the two It\^o projections give rise to new, distinct, Gaussian approximations. 

All our calculations of the resulting filters for the cubic sensor $b(x,t)=x + \epsilon x^3$ are equal when $\epsilon=0$. This provides a basic
sanity check that our formulae correspond to the Kalman filter in
the case of a linear sensor. In general, if we know that the solution
lies in a particular manifold and we project onto that manifold, the
three projection methods will all be exact.

We simulated the example problem $b(x)=x + \epsilon x^3$ for all of the above
approximate filters with $\epsilon=0.05$. We also computed an ``exact'' solution
using a finite difference method on a grid of 1000 intervals spaced evenly from $-10.0$ to $10.0$ and a time step of $0.0002$. We define the $L^2$ residual
to be the $L^2$ distance between the approximate solution and the ``exact'' solution.
We define the Hellinger residual similarly, as the $L^2$ distance between the square roots of the solution densities.

\begin{figure}[tbph]
\centering
\includegraphics[width=\linewidth]{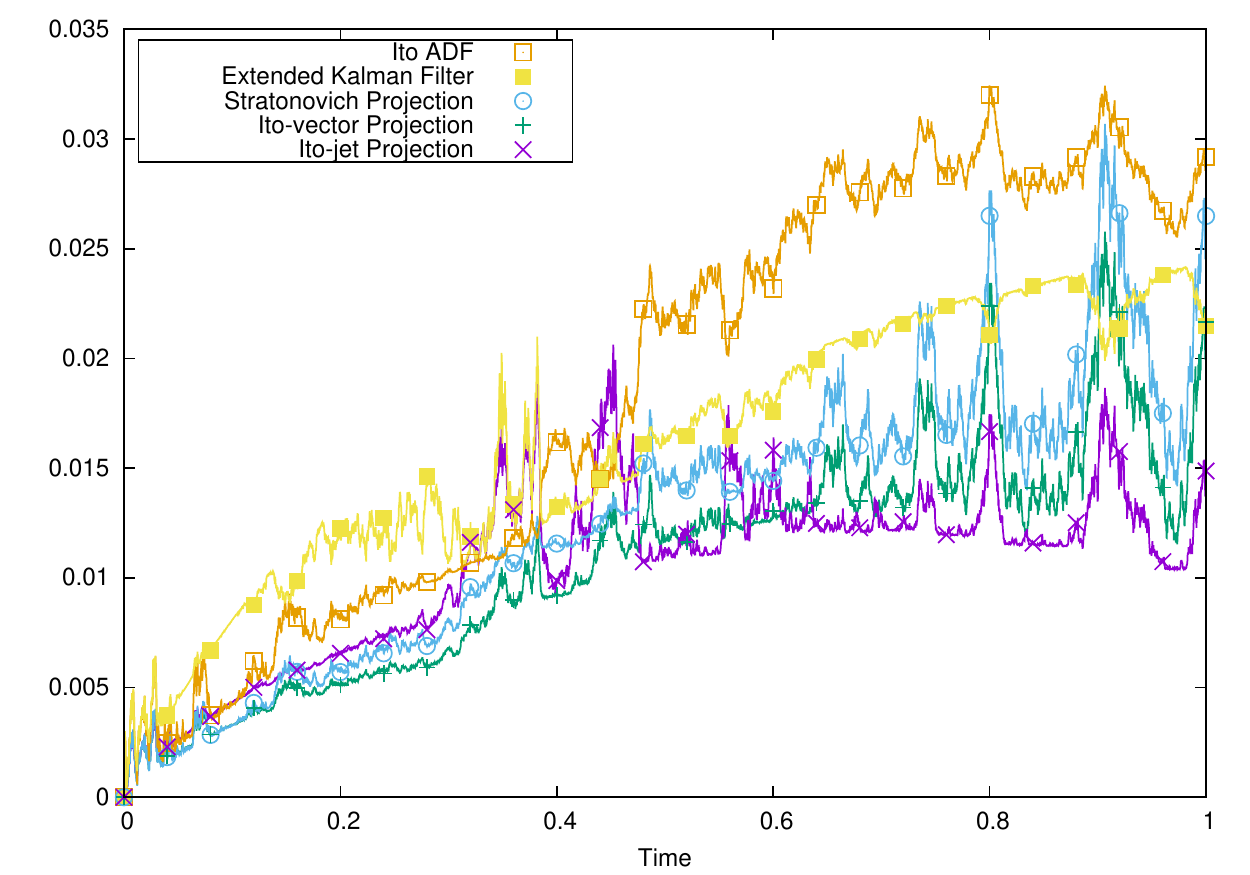}
\caption{$L^2$ residuals for each approximation method. All projections
are taken relative to the $L^2$ metric.}
\label{l2Residuals}
\end{figure}

In Figure \ref{l2Residuals} we see the $L^2$ residuals for the various methods.
All the projection methods shown are taken using the $L^2$ metric in this case.
The It\^o-vector projection in the $L^2$ metric results in the lowest
residuals over short time horizons. The Stratonovich projection comes
a close second. Over medium term time horizons, the It\^o-jet projection
out performs the It\^o-vector projection. We have not shown longer term behaviour because over long time horizons, all the methods become inaccurate and any comparison becomes meaningless. The projection methods out-performed all
other methods.
Although our plot shows only a single run, it is reasonably representative of the typical behaviour.

\begin{figure}[t]
\centering
\includegraphics[width=\linewidth]{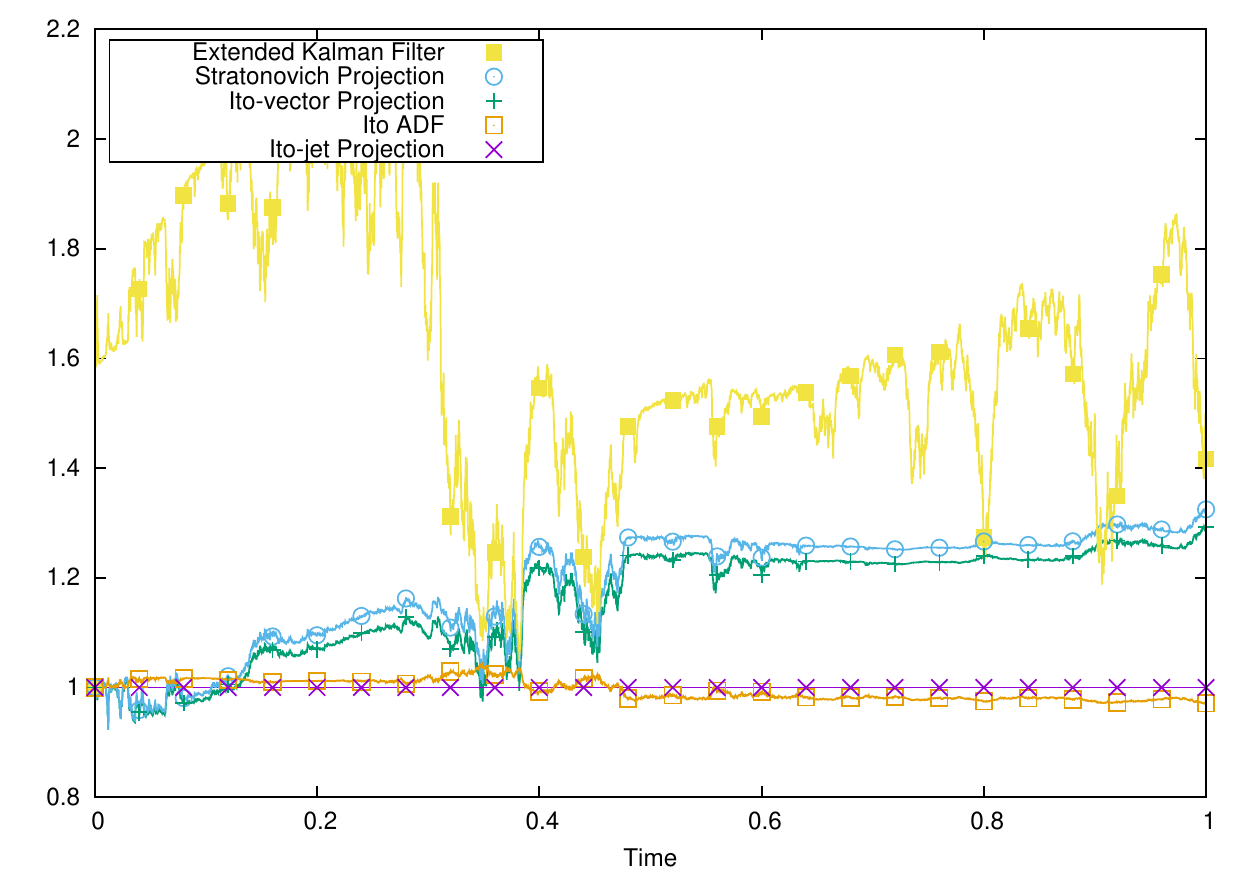}
\caption{Hellinger residuals for various approximation method divided by residual for the It\^o-jet projection. All projections
are taken relative to the Hellinger metric.}
\label{hellingerResiduals}
\end{figure}

In \Cref{hellingerResiduals} we have plotted the ratio of the Hellinger residual for
each method to the residual of the It\^o-jet projection w.r.t.
the Hellinger metric. This is because
the residuals themselves are too difficult to distinguish visually.
Thus values exceeding $1$ show a larger error than the It\^o-jet
projection and values less than one show a lower error. All the projection
methods shown in this plot are taken w.r.t.\ the Hellinger metric.

This plot indicates that the It\^o ADF and the It\^o-jet projection
are almost indistinguishable in their performance. A look at the
explicit formulae reveals that the difference between these two equations
is of order $\epsilon^2$ whereas the difference between the other equations
is of order only $\epsilon$.
Over the short term, the It\^o-vector projection gives the best
results. Over the medium term, the
It\^o-jet projection and the It\^o ADF give the best results.
Again, over the longer term all the filters become highly inaccurate.

\section{Conclusions}
\label{section:conclusions}

The notion of projecting a vector field onto a manifold is unambiguous.
By contrast, there are multiple distinct generalizations of this notion to
SDEs, as summarized in \Cref{tab:projtypes}.

The two It\^o projections we introduced in this work can both be derived
from minimization arguments. However, the It\^o-jet projection has some clear advantages.
\begin{itemize}
\item  The It\^o-jet projection is the best approximation to the metric projection of the true solution and has an error of $O(t)$. By contrast the It\^o-vector projection only tracks the true solution an accuracy of $O(t^{\frac{1}{2}})$.
\item  The It\^o-jet projection gives a more intuitive answer than
	   the It\^o-vector projection for the
	   low dimensional example considered in \Cref{section:lowDimensionalExample}.
\item  The It\^o-jet projection gives better numerical results in the
	   medium term than the It\^o-vector projection in our application to filtering.	   
\item  The It\^o-jet projection has an elegant definition when written in terms
       of $2$-jets.
\item  The It\^o-jet projection has a pictorial interpretation, shown in \Cref{fig:projectedSDE}.       	   
\end{itemize}

We have also seen that the Stratonovich projection satisfies an ad hoc minimization that is less appealing than the ones of the It\^o projections, since it requires a deterministic anchor point. 
%
%
%
The It\^o-jet and  It\^o-vector projection arguments allow one
to derive new Gaussian approximations to non-linear filters.
Unlike
previous Gaussian approximations to non-linear filters, these approximations
are derived by minimization arguments rather than heuristic arguments.
Thus the notion of projecting an SDE onto a manifold is able to give
new results even for this well-worn topic.

\appendix

\section{Appendix: The cross-diffusion process}\label{appendixcross}

We briefly study and give some intuition on the cross-diffusion process \eqref{crossDiffusion}, whose equations we repeat here:
\begin{equation*}
\begin{split}
\ed X_t &= \sigma Y_t \ dW_t, \\
\ed Y_t &= \sigma X_t \ dW_t,
\end{split}
\end{equation*}
with deterministic initial condition $(X_0,Y_0)$. We call this a cross diffusion, since each state crosses over as diffusion coefficient of the other state. Moreover, as we explain below, depending on the location on the plane of the initial condition, the paths group around the left or right arms of a St Andrew cross. This is another reason for the name of the process.

\begin{figure}[h!]
\centering
\includegraphics[width=0.8\linewidth]{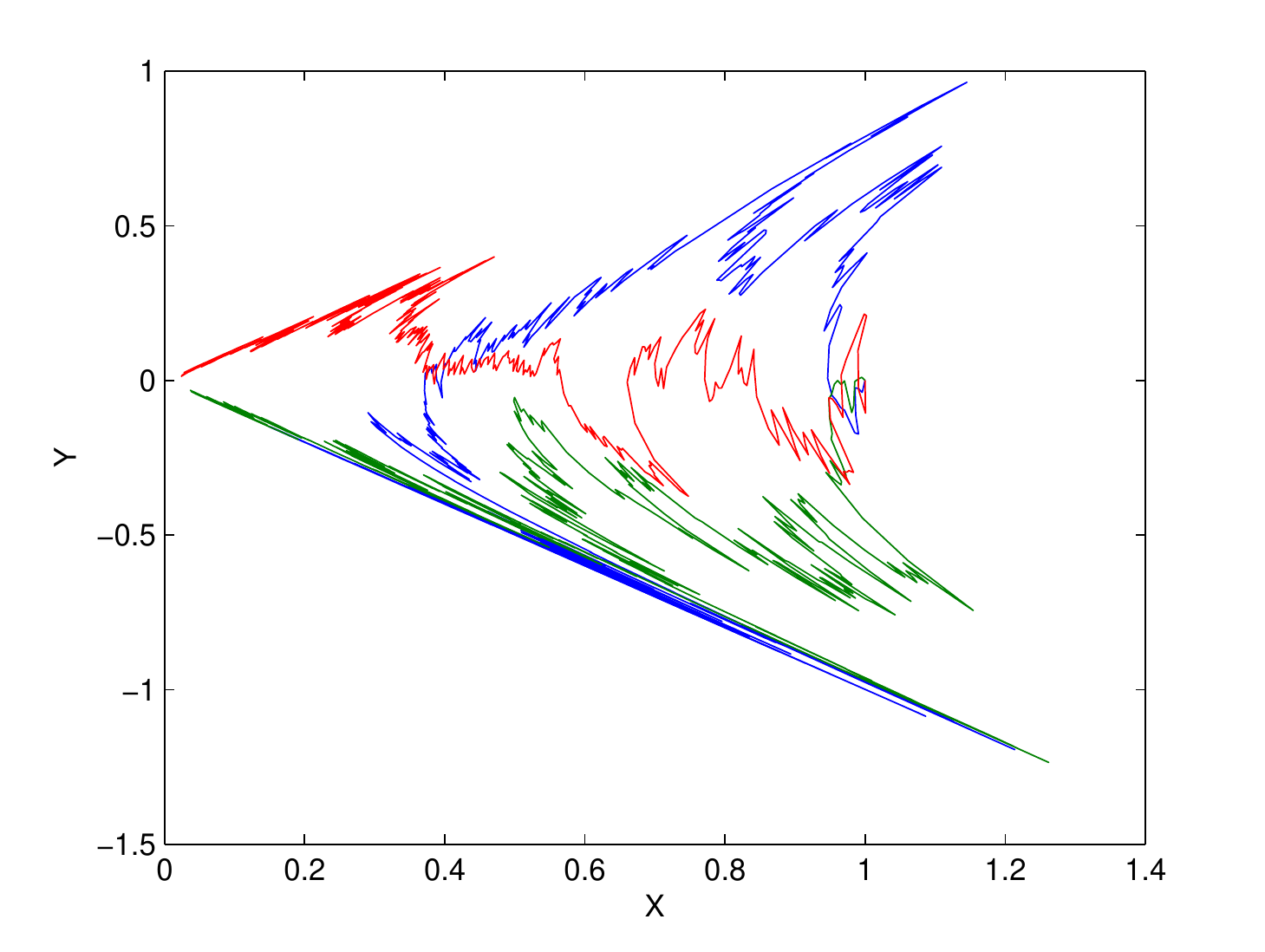}
\includegraphics[width=0.8\linewidth]{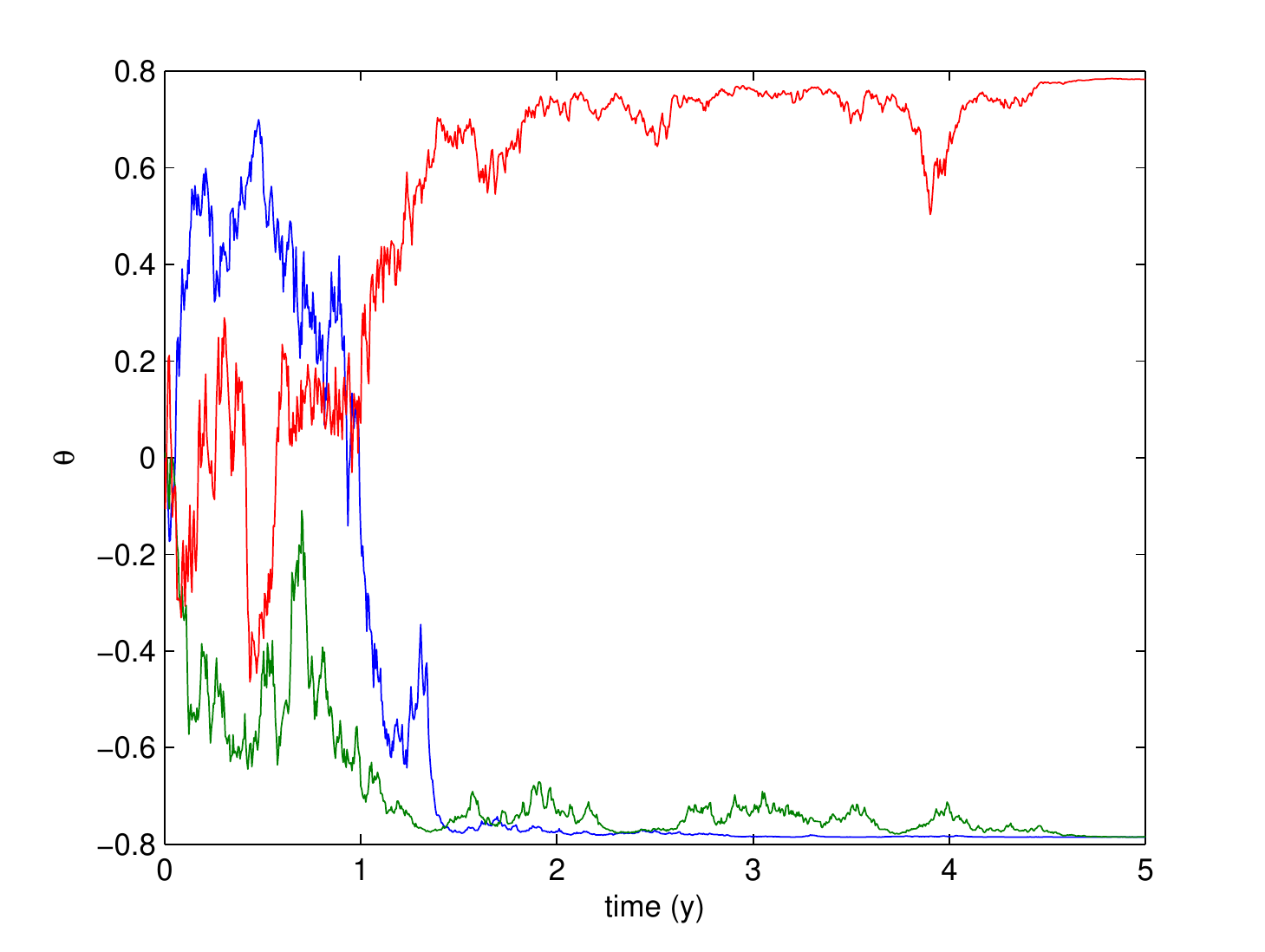}
\caption{Top: Three paths of the cross-diffusion SDE with $(X_0,Y_0) = (1,0)$, $\sigma=1$, up to $5$ years time. Bottom:  The corresponding three paths for $\theta$ plotted against time. }
\label{fig:crosspaths}
\end{figure}

The process equation can be solved analytically: add and subtract both sides and solve the resulting geometric Brownian motion equations for $X+Y$ and $X-Y$. One obtains:
\[
\begin{split}
X_t &= e^{-\frac{1}{2}\sigma^2 t} \left(X_0 \cosh(\sigma W_t) + Y_0 \sinh( \sigma W_t) \right), \\
Y_t &= e^{-\frac{1}{2}\sigma^2 t} \left(Y_0 \cosh(\sigma W_t) + X_0 \sinh( \sigma W_t )\right).
\end{split}
\]
One can see that the solution satisfies $(X_t+Y_t)(X_t-Y_t) = K e^{-\sigma^2 t}$, with $K=(X_0+Y_0)(X_0-Y_0)$. For large times the product will tend to be closer and closer to zero, so that either $Y=X$ or $Y=-X$, with the solution paths grouping along these two lines while approaching zero. Notice that if (X,Y) is in the origin, the process does not move. This behaviour of the process can be easily seen when plotting a few paths. In  \Cref{fig:crosspaths} we show a few paths of the cross-diffusion example in coordinates $(x,y)$ and the process $\theta_t = \arctan(Y_t/X_t)$ for the angular position. Clearly, when $X+Y$ is near zero $\theta$ will tend to be close to $-\pi/4 = - 0.785\ldots$, whereas when $X-Y$ is near zero $\theta$ will tend to be close to $\pi/4 = 0.785\ldots$.

\newpage

\bibliographystyle{plain}
\bibliography{filtering}

\begin{thebibliography}{10}

\bibitem{ahmedbook}
N.~U. Ahmed.
\newblock {\em Linear and Nonlinear Filtering for Scientists and Engineers}.
\newblock World Scientific, Singapore, 1998.

\bibitem{armstrongBrigoJets}
J.~Armstrong and D.~Brigo.
\newblock Coordinate free stochastic differential equations as jets.
\newblock {\em http://arxiv.org/abs/1602.03931}, 2016.

\bibitem{armstrongBrigo}
John Armstrong and Damiano Brigo.
\newblock Stochastic filtering via {L}2 projection on mixture manifolds with
  computer algorithms and numerical examples.
\newblock {\em arXiv preprint arXiv:1303.6236}, 2013.

\bibitem{armstrongBrigoicms}
John Armstrong and Damiano Brigo.
\newblock Extrinsic projection of {I}t{\^o} {SDE}s on submanifolds with
  applications to non-linear filtering.
\newblock {\em In: Nielsen, F., Critchley, F., \& Dodson, K. (Eds),
  Computational Information Geometry for Image and Signal Processing, Springer
  Verlag}, 2016.

\bibitem{armstrongbrigomcss}
John Armstrong and Damiano Brigo.
\newblock {N}onlinear filtering via stochastic {PDE} projection on mixture
  manifolds in {$L^2$} direct metric.
\newblock {\em Mathematics of Control, Signals and Systems}, 28(1):1--33, 2016.

\bibitem{crisan}
Alan Bain and Dan Crisan.
\newblock {\em Fundamentals of stochastic filtering}, volume~3.
\newblock Springer, 2009.

\bibitem{brewin}
L.~C. Brewin.
\newblock Riemann normal coordinates.
\newblock {\em Preprint Department of Mathematics, Monash University, Clayton,
  Victoria}, 3168, 1997.

\bibitem{brigosmall}
Damiano Brigo.
\newblock On the nice behaviour of the gaussian projection filter with small
  observation noise.
\newblock {\em Syst. Control Lett.}, 26(5):363--370, 1995.

\bibitem{brigo1}
Damiano Brigo, Bernard Hanzon, and Fran{\c{c}}ois LeGland.
\newblock A differential geometric approach to nonlinear filtering: the
  projection filter.
\newblock {\em Automatic Control, IEEE Transactions on}, 43(2):247--252, 1998.

\bibitem{brigo2}
Damiano Brigo, Bernard Hanzon, and Fran{\c{c}}ois LeGland.
\newblock Approximate nonlinear filtering by projection on exponential
  manifolds of densities.
\newblock {\em Bernoulli}, 5(3):495--534, 1999.

\bibitem{brigopistone}
Damiano Brigo and Giovanni Pistone.
\newblock Dimensionality reduction for measure-valued evolution equations in
  statistical manifolds.
\newblock In {\em Proceedings of the conference on Computational Information
  Geometry for Image and Signal Processing}. Springer, 2016.

\bibitem{brigopistone2}
Damiano Brigo and Giovanni Pistone.
\newblock Optimal approximations of the {F}okker-{P}lanck-{K}olmogorov
  equation: projection, maximum likelihood eigenfunctions and {G}alerkin
  methods.
\newblock {\em arXiv preprint arXiv:1603.04348}, 2016.

\bibitem{brzezniak}
Z.~Brze\'zniak and K.~D. Elworthy.
\newblock Stochastic differential equations on {B}anach manifolds.
\newblock {\em Methods Funct. Anal. Topology}, 6(1):43--84, 2000.

\bibitem{elworthy}
David Elworthy.
\newblock Geometric aspects of diffusions on manifolds.
\newblock In {\em {\'E}cole d'{\'E}t{\'e} de Probabilit{\'e}s de Saint-Flour
  XV--XVII, 1985--87}, pages 277--425. Springer, 1988.

\bibitem{emery}
M.~Emery.
\newblock {\em Stochastic calculus in manifolds}.
\newblock Springer-Verlag, Heidelberg, 1989.

\bibitem{girolami}
Mark Girolami and Ben Calderhead.
\newblock Riemann manifold {L}angevin and {H}amiltonian {M}onte {C}arlo
  methods.
\newblock {\em Journal of the Royal Statistical Society: Series B (Statistical
  Methodology)}, 73(2):123--214, 2011.

\bibitem{gliklikh}
Yuri~E. Gliklikh.
\newblock {\em {Global and Stochastic Analysis with Applications to
  Mathematical Physics}}.
\newblock Theoretical and Mathematical Physics. Springer, London, 2011.

\bibitem{hazewinkel}
Michiel Hazewinkel, SI~Marcus, and HJ~Sussmann.
\newblock Nonexistence of finite-dimensional filters for conditional statistics
  of the cubic sensor problem.
\newblock {\em Systems \& control letters}, 3(6):331--340, 1983.

\bibitem{hsu}
Elton~P. Hsu.
\newblock {\em Stochastic Analysis on Manifolds}.
\newblock Contemporary Mathematics. American Mathematical Society, 2002.

\bibitem{ito2}
K.~It{\^o}.
\newblock {S}tochastic differential equations in a differentiable manifold.
\newblock {\em Nagoya Math. J.}, (1):35--47, 1950.

\bibitem{jazwinski}
A.~H. Jazwinski.
\newblock {\em Stochastic Processes and Filtering Theory}.
\newblock Academic Press, New York, 1970.

\bibitem{kloedenAndPlaten}
Peter~E. Kloeden and Eckhard Platen.
\newblock {\em Numerical solution of stochastic differential equations}.
\newblock Applications of mathematics. Springer, Berlin, New York, Third
  printing, 1999.

\bibitem{kushner}
Harold~J Kushner.
\newblock Approximations to optimal nonlinear filters.
\newblock {\em Automatic Control, IEEE Transactions on}, 12(5):546--556, 1967.

\bibitem{maybeck}
Peter~S Maybeck.
\newblock {\em Stochastic models, estimation, and control}, volume~3.
\newblock Academic press, 1982.

\bibitem{newton1}
Nigel~J. Newton.
\newblock An infinite-dimensional statistical manifold modelled on {H}ilbert
  space.
\newblock {\em {J}. {F}unct. {A}nal.}, 263(6):1661--1681, 2012.

\bibitem{newton2}
Nigel~J. Newton.
\newblock Information geometric nonlinear filtering.
\newblock {\em {I}nfin. {D}imens. {A}nal. {Q}uantum {P}robab. {R}elat. {T}op.},
  2(18):1550014, 24, 2015.

\bibitem{pardoux}
{\'E}tienne Pardoux.
\newblock Filtrage non lin{\'e}aire et {\'e}quations aux d{\'e}riv{\'e}es
  partielles stochastiques associ{\'e}es.
\newblock In {\em Ecole d'Et{\'e} de Probabilit{\'e}s de Saint-Flour XIX 1989},
  pages 68--163. Springer, 1991.

\bibitem{picard}
Jean Picard.
\newblock Efficiency of the extended kalman filter for nonlinear systems with
  small noise.
\newblock {\em SIAM Journal on Applied Mathematics}, 51(3):843--885, 1991.

\bibitem{pistoneannals}
Giovanni Pistone and Carlo Sempi.
\newblock An infinite-dimensional geometric structure on the space of all the
  probability measures equivalent to a given one.
\newblock {\em Ann. Statist.}, 23(5):1543--1561, 10 1995.

\end{thebibliography}

\end{document}